\documentclass[12pt,a4paper]{article}

\usepackage{mathtools}
\usepackage{graphicx}   %
\usepackage{enumerate}  %
\usepackage{latexsym}   %
\usepackage{amsmath}    %
\usepackage{amsbsy}     %
\usepackage{amsfonts}   %
\usepackage{mathrsfs}   %
\usepackage{amssymb}    %
\usepackage{array}      %
\usepackage{theorem}    %
\usepackage{mathrsfs,bm}
\usepackage{comment}
\usepackage{color}

\usepackage[textsize=small]{todonotes}
 \newtheorem{thm}{Theorem}[section]
 
 \newtheorem{lem}[thm]{Lemma}
 
 \theoremstyle{definition}
 \newtheorem{df}[thm]{Definition}
 \theoremstyle{remark}
 \newtheorem{rem}[thm]{Remark}
 
 \numberwithin{equation}{section}
 \newtheorem{notn}{Notation}

\def\D{{\mathscr D}}
\def\F{{\mathscr F}}
\def\N{{\mathbb N}}
\def\R{{\mathbb R}}

\def\J{{\mathscr J}}

\def\v{\bm{v}}
\def\vv{\bm{v}^{\nu,\sigma}}
\def\thetaa{\theta^{\nu,\sigma}}
\def\w{\bm{w}}
\def\ww{w^{\nu,\sigma}}
\def\ds{\displaystyle}
\def\ts{\textstyle}

\def\laplace{\varDelta}

\def\laplace{{\varDelta}}

\definecolor{wineRed}{rgb}{0.7,0,0.3}
\definecolor{grandBleu}{rgb}{0,0,0.8}
\definecolor{darkGreen}{rgb}{0,0.5,0}
\definecolor{blueViolet}{rgb}{0.4,0,1.0}
\definecolor{bloodOrange}{rgb}{0.85,0.05,0}
\definecolor{mycolor}{rgb}{0.8,0,0.2}
\definecolor{magenta}{rgb}{0.8,0,1.0}
\definecolor{}{rgb}{0.8,0,0.2}

\begin{document}
\label{page:t}

\vspace*{1cm}
\begin{center}
  {\large \bf {\bf
          Energy-dissipation in a coupled system of Allen--Cahn  \\ type equation and Kobayashi--Warren--Carter \\[0.5ex] type model of grain boundary motion\footnotemark[1]
  }}
\end{center}
\vspace{7ex}
\begin{center}

{\sc Hiroshi Watanabe}\footnotemark[2]\\
Division of Mathematical Sciences,\\
Faculty of Science and Technology, Oita University\\
700 Dannoharu, Oita, 870--1192, Japan\\
{\ttfamily hwatanabe@oita-u.ac.jp}
\vspace{7ex}

{\sc Ken Shirakawa}\footnotemark[3]\\
Department of Mathematics, Faculty of Education, Chiba University\\
1--33, Yayoi-cho, Inage-ku, Chiba, 263--8522, Japan\\
{\ttfamily sirakawa@faculty.chiba-u.jp}

\end{center}
\vspace{14ex}


\noindent
{\bf Abstract. } 
In this paper, we consider a system of initial boundary value problems for parabolic equations, as a generalized version of the ``$ \phi $-$ \eta $-$ \theta $ model'' of grain boundary motion, proposed by Kobayashi \cite{K}. The system is a coupled system of: an Allen--Cahn type equation as in \eqref{eqn-2} with a given temperature source; and a phase-field model of grain boundary motion, known as ``Kobayashi--Warren--Carter type model''. The focus of the study is on a special kind of solution, called {\em energy-dissipative solution}, which is to reproduce the energy-dissipation of the governing energy in time. Under suitable assumptions, two Main Theorems, concerned with: the existence of energy-dissipative solution; and the large-time behavior; will be demonstrated as the results of this paper. 


\footnotetext[0]{
$\empty^*$\,AMS Subject Classification 35K87, 
	35R06, 
	35K67.  

    $\empty^*$\,Keywords: Allen--Cahn type equation; Kobayashi--Warren--Carter type model; grain boundary motion; smoothing effect; energy dissipation; large-time behavior.

$\empty^\dag$\,This author is supported by Grant-in-Aid No. 16K05224, JSPS.
$\empty^\ddag$\,This author is supported by Grant-in-Aid No. 17K05294, JSPS. 
}
\newpage
\section{Introduction}
\ \ \vspace{-4ex}

Let $ N \in \N $ be a constant of spatial dimension, and $ \Omega \subset \R^N $ be a bounded domain such that $ \Gamma := \partial \Omega $ is smooth when $ N > 1 $. Besides, let us denote by $ Q :=  (0, \infty) \times \Omega $ the product space of the time-interval $ (0, \infty) $ and the spatial domain $ \Omega $, and similarly, let us set $ \Sigma := (0, \infty) \times \Gamma $.  

In this paper, we fix a constant $ \nu \geq 0 $, and consider the following system of initial-boundary value problems of parabolic types, denoted by (S)$_\nu$.
\bigskip

\noindent 
(S)$_\nu$:
\begin{align}
    \begin{cases}
        \parbox{11.5cm}{
            $ \ds w_{t} - \laplace w + \partial \gamma(w) +g_w(w, \eta) +cu $ 
            \\
            \hspace*{8.7ex} $ + \alpha_{w}(w,\eta) |D\theta| + \nu^2 \beta_{w}(w,\eta)|D\theta|^{2} \ni 0 \mbox{ \ in } Q, $ 
            \\[1.5ex]
            $ \ds D w \cdot \bm{n}_\Gamma = 0\ \ \mbox{ on } \Sigma $,
            \\[1.5ex]
            $ w(0, x)= w_{0}(x), ~~ x \in \Omega $;
        }
    \end{cases}
    \label{eqn-2}
\end{align}
\begin{align}
    \begin{cases}
        \parbox{11.5cm}{
            $ \ds \eta_{t} - \laplace \eta + g_{\eta} (w,\eta) + \alpha_{\eta}(w,\eta)|D\theta| +\nu^2 \beta_{\eta}(w,\eta)|D\theta|^{2} = 0\ \ \mbox{ in } Q $, 
            \\[1.5ex]
            $ \ds D \eta \cdot \bm{n}_\Gamma = 0\ \ \mbox{ on } \Sigma $,
            \\[1.5ex]
            $ \ds \eta(0, x)= \eta_{0}(x), ~~ x \in \Omega $; 
        }
    \end{cases}
    \label{eqn-3}
\end{align}
\begin{align}
    \begin{cases}
        \parbox{11.5cm}{
            $ \ds \alpha_{0}(w,\eta) \theta_{t} - \mbox{div} \left( \alpha(w,\eta) \frac{D\theta}{|D\theta|} + 2 \nu^2 \beta(w,\eta)D\theta \right) = 0\ \ \mbox{ in } Q $, 
            \\[1.5ex]
            $ \ds \bigl( {\textstyle \alpha(w,\eta) \frac{D\theta}{|D\theta|} +2 \nu^2 \beta(w,\eta)D\theta }\bigr) \cdot \bm{n}_\Gamma = 0\ \ \mbox{ on } \Sigma $, 
            \\[1.5ex]
            $ \ds \theta(0, x) = \theta_{0}(x), ~~ x \in \Omega $. 
        }
    \end{cases}
    \label{eqn-4}
\end{align}

The system (S)$_\nu$ is a generalized version of the ``$ \phi $-$ \eta $-$ \theta $ model'' of grain boundary motion, which was proposed by Kobayashi \cite{K}. The first initial-boundary value problem \eqref{eqn-2} is a type of Allen--Cahn equation, i.e. \eqref{eqn-2} is a mathematical model of solid-liquid phase transition in a polycrystal. Meanwhile, the system of second-third problems \{ \eqref{eqn-3}, \eqref{eqn-4} \} forms a type of Kobayashi--Warren--Carter model of grain boundary motion, which is proposed in \cite{KWC1,KWC2}, and studied by a lot of mathematicians from various viewpoints (e.g., \cite{GG,IKY08,IKY09,IKY12,KY,KG,KWC1,KWC2,MS14,MSW,SW14,SW15,SW17,SWY13,SWY14,SWY16,WKLC,WS13,WS14}).
    
The system (S)$_{\nu}$ is derived as a gradient system of the following governing energy, called ``{\em free-energy}":
\begin{equation}\label{freeEnrgy}
\begin{array}{ll}
	\multicolumn{2}{l}{\ds \mathscr{E}_\nu^{u}(w, \eta, \theta) := \frac{1}{2} \int_\Omega |D w|^2 \, dx +\int_\Omega \gamma(w) \, dx + c\int_\Omega u w \, dx}
\\[2.5ex]
& \ds +\frac{1}{2} \int_\Omega |D \eta|^2 \, dx +\int_\Omega g(w, \eta) \, dx  +\int_\Omega \alpha(w, \eta) \, d|D \theta| +\int_\Omega \beta(w, \eta)|D(\nu \theta)|^2 \, dx,
\\[2.5ex]
& \mbox{for $ [w, \eta, \theta] \in H^1(\Omega) \times H^1(\Omega) \times BV(\Omega) $ with $ \nu \theta \in H^1(\Omega) $.}
\end{array}
\end{equation}
In this context, $ u = u(t, x) $ is a given temperature source (relative temperature), and the unknown $ w = w(t, x) $ is an order parameter to indicate the solidification order of the polycrystal. The unknowns $ \eta = \eta(t, x) $ and $ \theta = \theta(t, x) $ are components of the vector field
$$
(t, x) \in Q \mapsto \eta(t, x) \left[ \rule{0pt}{10pt} \cos \theta(t, x), \sin \theta(t, x) \right] \in \mathbb{R}^2,
$$
which was adopted in \cite{KWC1,KWC2} as a vectorial phase-field to reproduce the crystalline orientation in $ Q $. Besides, the components $ \eta $ and $ \theta $ are order parameters to indicate, respectively, the orientation order and orientation angle of the grain. In particular, $ w $ and $ \eta $ are taken to satisfy the constraints $ 0 \leq w, \eta \leq 1 $ in $ Q $, and the cases $ [w, \eta] \approx [1, 1] $ and $ [w, \eta] \approx [0, 0] $ are respectively assigned to ``the solidified-oriented phase'' and ``the liquefied-disoriented phase'' which correspond to two stable phases in physics. Meanwhile, $ w_0 = w_0(x) $, $ \eta_0 = \eta_0(x)$ and $ \theta_0 = \theta_0(x) $ are given initial data on $ \Omega $. $ \partial \gamma $ is the subdifferential of a proper lower semi-continuous (l.s.c.) and convex function $ \gamma = \gamma(w) $ on $ \R $. $ u = u(x,t) $, $ g = g(w, \eta) $, $ \alpha_0 = \alpha_0(w, \eta) $, $ \alpha = \alpha(w, \eta) $, and $ \beta = \beta(w, \eta) $ are given real-valued functions, and the scripts ``$ \empty_w $'' and ``$ \empty_\eta $'' denote differentials with respect to the corresponding variables. $ \bm{n}_\Gamma $ is the unit outer normal on $ \Gamma $.

With regard to the Kobayashi--Warren--Carter type models, the most of mathematical results, obtained in the previous works \cite{MS14,MSW,SW14,SW15,SW17,SWY13,SWY14,SWY16,WS13,WS14}, are classified in the following four issues. 
\begin{description}
\item[(T1)]Variational solvability, i.e. the existence of solution in the variational sense.
\item[(T2)]Existence of ``weak solution'', which is to realize the ``smoothing effect'' as a solution to a parabolic system.
\item[(T3)]Existence of ``energy-dissipative solution'', which is to realize the ``energy-dissipation'', i.e. the nonincreasing property associated with the time-variation of free-energy.
\item[(T4)]Large time behavior of the energy-dissipative solution. 
\end{description}

For the original Kobayashi--Warren--Carter model, the mathematical results concerned with (T1)--(T4) were studied in  \cite{MS14,MSW,SW14,SWY13,WS13,WS14} under suitable assumptions. Meanwhile, in the mathematical analysis for the system (S)$_\nu$, we still have some incomplete parts. More precisely, for a simplified version of (S)$_\nu$, the issues (T1), (T3) and (T4) were studied in \cite{SW15}, and the result is extended to the mathematical analysis under unknown setting of the temperature $ u $ (cf. \cite{SW17}). However, for general case of (S)$_\nu$, there is only one result for (T1) (cf. \cite{SWY14}), and there is no result to give mathematical answers for the remaining issues (T2)--(T4), yet. 

In view of such background, we set the goal of this paper to establish a general mathematical theory that enable a uniform treatment for the issues (T1)--(T4), under various settings of the system (S)$_\nu$. On this basis, the principal discussion  will be devoted to the proofs of the following two main theorems. 

\begin{description}
	\item[Main Theorem 1:]the existence theorem of energy-dissipative solutions $ [w, \eta, \theta] $ to the systems (S)$_\nu$, for any $ \nu \geq 0 $, which behaves in the range of $ C([0, \infty); L^2(\Omega)^3) $.
	\item[Main Theorem 2:]the large-time behavior of energy-dissipative solutions.
\end{description}

    The contents of this paper are as follows. The Main Theorems are stated in Section 3, after the preliminaries in Section 2. The Main Theorems are proved in the following Sections 5 and 6, and in particular, the proof of Main Theorem 1 is based on some Lemmas for approximation problem, obtained in Section 4.

\section{Preliminaries}
First we elaborate the notation used throughout. 
\begin{notn}[Abstract notations]\label{Note00}
\begin{em}
Let $ d \in \N $ take any fixed value. The $ d $-dimensional Lebesgue measure is denoted by $ \mathscr{L}^d $. Also, unless otherwise specified, the measure-theoretic phrases such as ``a.e.,'' ``$ dt $,'' ``$ dx $'', and so on, are with respect to the Lebesgue measure in each corresponding dimension.
\end{em}
\end{notn}

\begin{notn}[Abstract functional analysis]\label{Note01}
\begin{em}
For an abstract Banach space $ X $, we denote by $ |\mathrel{\cdot}|_X $ the norm of $ X $, and when $ X $ is a Hilbert space, we denote by $ (\,\cdot\,, \,\cdot\,)_X $ its inner product. For a subset $ A $ of a Banach space $ X $, we denote by $ {\rm int}(A) $ and $ \overline{A} $ the interior and the closure of $ A $, respectively. 

Fix $ 1 < d \in \N $. Then, for a Banach space $ X $ the topology of the product Banach space
\begin{center}
$ X^d := \overbrace{X \times \cdots \times X}^{\mbox{\scriptsize $ d $ times}} $
\end{center}
has the norm
\begin{equation*}
|z|_{X^d} := \sum_{k = 1}^d |z_k|_{X}, \mbox{ for $ z = [z_1, \ldots, z_d] \in X^d $.}
\end{equation*}
However, if $ X $ is a Hilbert space, then the topology of the product Hilbert space $ X^d $ has the inner product
\begin{equation*}
(z, \tilde{z})_{X^d} := \sum_{k = 1}^d (z_k, \tilde{z}_k)_X, \mbox{ for $ z = [z_1, \ldots, z_d] \in X^d $ and $ \tilde{z} = [\tilde{z}_1, \ldots, \tilde{z}_d] \in X^d $,}
\end{equation*}
and hence, the norm in this case is provided by
\begin{equation*}
|z|_X := \sqrt{(z, z)_{X^d}} = \left( \rule{-2pt}{18pt} \right. \sum_{k = 1}^d |z_k|_X^2 \left. \rule{-2pt}{18pt} \right)^{\hspace{-0.5ex}1/2}, \mbox{ for $ z = [z_1, \ldots, z_d] \in X^d $.}
\end{equation*}

For any proper lower semi-continuous (l.s.c.\ hereafter) and convex function $ \Psi $ defined on a Hilbert space $ X $, we denote by $ D(\Psi) $ its effective domain, and by $ \partial \Psi $ its subdifferential. The subdifferential $ \partial \Psi $ is a set-valued map corresponding to a weak differential of $ \Psi $, and it has a maximal monotone graph in the product Hilbert space $ X^2 $. More precisely, for each $ z_0 \in X $, the value $ \partial \Psi(z_0) $ is defined as the set of all elements $ z_0^* \in X $ that satisfy the variational inequality
\begin{equation*}
(z_0^*, z -z_0)_X \leq \Psi(z) -\Psi(z_0) \mbox{ \ for any $ z \in D(\Psi) $,}
\end{equation*}
and the set $ D(\partial \Psi) := \{z \in X \mid \partial \Psi(z) \ne \emptyset\} $ \
is called the domain of $ \partial \Psi $. We often use the notation ``$ [z_0, z_0^*] \in \partial \Psi $ in $ X^2 $\,'' to mean ``$ z_0^* \in \partial \Psi(z_0) $ in $ X $ with $ z_0 \in D(\partial \Psi) $,'' by identifying the operator $ \partial \Psi $ with its graph in $ X^2 $.
\end{em}
\end{notn}
\begin{rem}\label{Rem.Time-dep.}
\begin{em}
It is often useful to consider the subdifferentials under time-dependent settings. In this regard, several general theories have been established by previous researchers (e.g., Kenmochi \cite{Kenmochi01}, and \^{O}tani \cite{Otani}). From these (e.g., \cite[Chapter 2]{Kenmochi01}), one can see the following fact: 
\begin{description}
\item[{(Fact1)}]Let $ E_0 $ be a convex subset in a Hilbert space $ X $, let $ I \subset [0, \infty) $ be a time interval, and for any $ t \in I $, let $ \Psi^t : X \rightarrow (-\infty, \infty] $ be a proper l.s.c. and convex function such that $ D(\Psi^t) = E_0 $ for all $ t \in I $. Based on this, define a convex function $ {\Psi}^I : L^2(I; X) \rightarrow (-\infty, \infty] $, by setting
\begin{equation*}
\zeta \in L^2(I; X) \mapsto {\Psi}^I(\zeta) := \left\{ \begin{array}{ll}
\multicolumn{2}{l}{\ds \int_I \Psi^t(\zeta(t)) \, dt, \mbox{ if $ \Psi^{(\cdot)}(\zeta) \in L^1(I) $,}}
\\[1ex]
\infty, & \mbox{otherwise.}
\end{array} \right.
\end{equation*}
Here, if $ E_0 \subset D({\Psi}^I) $, and the function $ t \in I \mapsto \Psi^t(z) $ is integrable for any $ z \in E_0 $, then the following holds:
\begin{equation*}
\begin{array}{c}
[\zeta, \zeta^*] \in \partial {\Psi}^I \mbox{ \ in $ L^2(I; X)^2 $ if and only if}
\\[1ex]
\zeta \in D({\Psi}^I) \mbox{ and } [\zeta(t), \zeta^*(t)] \in \partial \Psi^t \mbox{ in $ X^2 $, a.e.\ $ t \in I $.}
\end{array}
\end{equation*}
\end{description}
\end{em}
\end{rem}
\begin{notn}[Basic elliptic operators]\label{Note02}
\begin{em}
Let $ {\mit \Delta}_{N} $ be the Laplacian operator subject to the zero Neumann boundary condition, i.e.,
\begin{equation*}
    {\mit \Delta}_{N} : z \in D(\mathit{\Delta}_N):= \left\{ \begin{array}{l|l}
z \in H^2(\Omega) & \nabla z \cdot \nu_{\partial \Omega} = 0 \mbox{ in $ L^2(\partial \Omega) $}
\end{array} \right\} \subset L^2(\Omega) \mapsto {\mit \Delta} z \in L^2(\Omega). 
\end{equation*}
Let $ d \in \N $ be a fixed constant of dimension. Then, we let:
\begin{equation*}
    \mathit{\Delta}_N z := \bigl[ \mathit{\Delta}_N z_1, \dots, \mathit{\Delta}_N z_d \bigr] \in L^2(\Omega)^d, \mbox{ \ for all $ z = [z_1, \dots, z_d] \in D(\mathit{\Delta}_N)^d $.}
\end{equation*}
As is well-known (see, e.g. \cite{Barbu} or \cite{Brezis}), the operator $ -\mathit{\Delta}_N $ coincides with the subdifferential of a proper l.s.c. and convex function  $ V_{D}^d : L^2(\Omega)^d \longrightarrow [0, \infty] $, defined as:
\begin{equation*}\label{Dirichlet^d}
z \in L^2(\Omega)^d \mapsto V_{D}^d(z) := \left\{ \begin{array}{ll}
\multicolumn{2}{l}{\ds \frac{1}{2} \int_\Omega |\nabla z|_{\R^{d \times N}}^2 \, dx, \mbox{ \ if $ z \in H^1(\Omega)^d $,}}
\\[2ex]
\infty, & \mbox{otherwise.}
\end{array} \right.
\end{equation*}
More precisely, 
\begin{equation*}\label{elliptic01}
z \in L^2(\Omega)^d \mapsto \partial V_{D}^d(z) = \left\{ \begin{array}{ll}
\multicolumn{2}{l}{\{ -{\mit \Delta}_{N} z \}, \mbox{ if $ z \in D_{N}^d $,}}
\\[1ex]
\emptyset, & \mbox{otherwise.}
\end{array} \right.
\end{equation*}
In this light, $ \partial V_{D}^d $ and $ -{\mit \Delta}_{N} $ are identified as the maximal monotone graphs in $ [L^2(\Omega)^d]^2 $. 
\end{em}
\end{notn}

\begin{notn}[BV theory; cf.\ \cite{AFP,ABM,EG,G}]\label{Note04}
\begin{em}
Let $ d \in \N $ be a fixed number, and let $U\subset \R^d$ be an open set. We denote by $ \mathcal{M}(U) $ (resp. $ \mathcal{M}_{\rm loc}(U) $) the space of all finite Radon measures (resp. the space of all Radon measures) on $ U $. In general, the space $ \mathcal{M}(U) $ (resp. $ \mathcal{M}_{\rm loc}(U) $) is known as the dual of the Banach space $ C_0(U) $ (resp. dual of the locally convex space $ C_{\rm c}(U) $), for any open set $ U \subset \mathbb{R}^d $.

A function $ v \in L^1(U) $ (resp. $ v \in L_{\rm loc}^1(U) $)  is called a function of bounded variation, or a BV-function, (resp. a function of locally bounded variation or a BV$\empty_{\rm loc}$-function) on $ U $, if and only if its distributional differential $ D v $ is a finite Radon measure on $ U $ (resp. a Radon measure on $ U $), namely $ D v \in \mathcal{M}(U) $ (resp. $ D v \in \mathcal{M}_{\rm loc}(U) $).
We denote by $ BV(U) $ (resp. $ BV_{\rm loc}(U) $) the space of all BV-functions (resp. all BV$\empty_{\rm loc}$-functions) on $ U $. For any $ v \in BV(U) $, the Radon measure $ D v $ is called the variation measure of $ v $, and its  total variation $ |Dv| $ is called the total variation measure of $ v $. Additionally, the value $|Dv|(U)$, for any $v \in BV(U)$, can be calculated as follows:
\begin{equation*}
|Dv|(U) = \sup \left\{ \begin{array}{l|l}
\ds \int_{U} v \ {\rm div} \,\varphi \, dy & \varphi \in C_{\rm c}^{1}(U)^d \ \ \mbox{and}\ \ |\varphi| \le 1\ \mbox{on}\ U
\end{array}
\right\}.
\end{equation*}
The space $BV(U)$ is a Banach space, endowed with the following norm:
\begin{equation*}
|v|_{BV(U)} := |v|_{L^{1}(U)} + |D v|(U),\ \ \mbox{for any}\ v\in BV(U).
\end{equation*}
We say that a sequence $\{ v_{n} \}_{n = 1}^\infty \subset BV(U)$ strictly converges in $BV(U)$ to $v \in BV(U)$ if $v_{n} \to v$ in $L^{1}(U)$ and $|Dv_{n}|(U) \to |Dv|(U)$ as $ n \to \infty $. Also, we say that a sequence $\{ v_{n} \}_{n = 1}^\infty \subset BV(U)$ (resp. $ \{ v_{n} \}_{n = 1}^\infty \subset BV_{\rm loc}(U) $) converges to $v \in BV(U)$ (resp. $ v \in BV_{\rm loc}(U) $) weakly-$*$ in $ BV(U) $ (resp. weakly-$*$ in $ BV_{\rm loc}(U) $) if $v_{n} \to v$ in $L^{1}(U)$ (in $ L_{\rm loc}^1(U) $) and $ Dv_{n} \to Dv $ weakly-$*$ in $ \mathcal{M}(U)^d $ (weakly-$*$ in $ \mathcal{M}_{\rm loc}(U)^d $) as $ n \to \infty $.
In particular, if the boundary $\partial U$ is Lipschitz, then the space $BV(U)$ is continuously embedded into $L^{d/(d-1)}(U)$ and compactly embedded into $L^{q}(U)$ for any $1 \le q < d/(d-1)$ (cf. \cite[Corollary 3.49]{AFP} or \cite[Theorem 10.1.3-10.1.4]{ABM}). Besides, any bounded subset in $ BV(U) $ (resp.  $ BV_{\rm loc}(U) $) is sequentially compact in $ BV(U) $ (resp.  $ BV_{\rm loc}(U) $) with respect to the weak-$*$ topology of $ BV(U) $ (resp. $ BV_{\rm loc}(U) $) (cf. \cite[Theorem 3.23]{AFP}). Additionally, if $1 \le r < \infty$, then the space $C^{\infty}(\overline{U})$ is dense in $BV(U) \cap L^{r}(U)$ for the intermediate convergence (cf. \cite[Definition 10.1.3. and Theorem 10.1.2]{ABM}), i.e. for any $v \in BV(U) \cap L^{r}(U)$, there exists a sequence $\{v_{n} \}_{n = 1}^\infty \subset C^{\infty}(\overline{U})$ such that $v_{n} \to v$ in $L^{r}(U)$ and $\int_{U}|\nabla v_{n}|dx \to |Dv|(U)$ as $n \to \infty$.
\end{em}
\end{notn}
\begin{notn}[Weighted total variation; cf.\ \cite{AB,AFP}]\label{Note05}
\begin{em}
In this paper, we define
\begin{equation*}
X_{\rm c}(\Omega) := \{
\varpi \in L^\infty(\Omega)^N  | \ \mbox{$ {\rm div} \, \varpi \in L^2(\Omega) $ \ and \ $ {\rm supp} \, \varpi $ is compact in $ \Omega $}
\},
\end{equation*}
\begin{equation*}
W_0(\Omega) := \left\{ \begin{array}{l|l}
\varrho \in H^1(\Omega) \cap L^\infty(\Omega) & \varrho \geq 0 \mbox{ a.e.\ in $ \Omega $}
\end{array} \right\},
\end{equation*}
\begin{equation}\label{W_c}
W_{\rm c}(\Omega) := \left\{ \begin{array}{l|l}
\varrho \in H^1(\Omega) \cap L^\infty(\Omega) & \parbox{5.5cm}{there exists $ c_\varrho > 0 $ such that $ \varrho \geq c_\varrho $ a.e.\ in $ \Omega $}
\end{array} \right\},
\end{equation}
and for any $ \varrho \in W_0(\Omega) $ and any $ z \in L^2(\Omega) $, we call the value $ {\rm Var}_\varrho(z) \in [0, \infty] $, defined as,
\begin{equation*}
{\rm Var}_\varrho(v) := \sup \left\{ \begin{array}{l|l}
\displaystyle \int_\Omega v \, {\rm div} \, \varpi \, dx & \parbox{3.5cm}{$ \varpi \in X_{\rm c} $, and \\ $ | \varpi | \leq \varrho $ a.e.\ in $ \Omega $}
\end{array} \right\} \in [0, \infty],
\end{equation*}
``the total variation of $ v $ weighted by $ \varrho $,'' or the ``weighted total variation'' for short. 
\end{em}
\end{notn}
\begin{rem}\label{Rem.Note05}
\begin{em}
Referring to the general theories (e.g., \cite{AB,AFP,BBF}), we can confirm the following facts associated with the weighted total variations: 
\begin{description}
\item[{(Fact2)}](cf.\ \cite[Theorem 5]{BBF}) For any $ \varrho \in W_0(\Omega) $, the functional $ z \in L^2(\Omega) \mapsto {\rm Var}_\varrho(z) \in [0, \infty] $ is a proper l.s.c. and convex function that coincides with the lower semi-continuous envelope of
\begin{equation*}
z \in W^{1, 1}(\Omega) \cap L^2(\Omega) \mapsto \int_\Omega \varrho |\nabla z| \, dx \in [0, \infty).
\end{equation*}
\item[{(Fact3)}](cf.\ \cite[Theorem 4.3]{AB} and \cite[Proposition 5.48]{AFP}) If $ \varrho \in W_0(\Omega) $ and $ z \in BV(\Omega) \cap L^2(\Omega) $, then there exists a Radon measure $ |Dz|_\varrho \in \mathcal{M}(\Omega) $ such that
\begin{equation*}
|Dz|_\varrho(\Omega) = \int_\Omega  d|Dz|_\varrho = {\rm Var}_\varrho(z),
\end{equation*}
and 
\begin{equation}\label{|Dz|_beta(A)}
|Dz|_\varrho(A) \leq |\varrho|_{L^\infty(\Omega)} |Dz|(A), \mbox{ for any open set }  A \subset \Omega.
\end{equation}
\item[{(Fact4)}]If $ \varrho \in W_{\rm c}(\Omega) $ and $ z \in BV(\Omega) \cap L^2(\Omega) $, then for any open set $ A \subset \Omega $, it follows that
\begin{equation}\label{|Dz|_BV}
\left\{ ~ {\parbox{10cm}{
$ |Dv|_\varrho(A) \geq c_\varrho |Dz|(A) $ \ for any open set $ A \subset \Omega $,
\\[2ex]
$ D({\rm Var}_\varrho) = BV(\Omega) \cap L^2(\Omega) $, and
\\[1ex]
$ {\rm Var}_\varrho(z) = \sup \left\{ \begin{array}{l|l}
\displaystyle \int_\Omega z \, {\rm div} \, (\varrho \varphi) \, dx & \parbox{3cm}{$ \varphi \in X_{\rm c}(\Omega) $, and $ | \varphi | \leq 1 $ a.e.\ in $ \Omega $}
\end{array} \right\}, $}} \right.
\end{equation}
where $ c_\varrho $ is a constant as in (\ref{W_c}). 
\end{description}
Moreover, the following properties can be inferred from (\ref{|Dz|_beta(A)})--(\ref{|Dz|_BV}):
\begin{itemize}
\item[$\bullet$] $ |Dz|_c = c|Dz| $ in $ \mathcal{M}(\Omega) $ for any constant $c \geq 0 $ and $ z \in BV(\Omega) \cap L^2(\Omega) $;
\item[$\bullet$] $ |Dz|_\varrho = \varrho |\nabla z| \mathscr{L}^N $ in $ \mathcal{M}(\Omega) $, \ if $ \varrho \in W_0(\Omega) $ and $ z \in W^{1, 1}(\Omega) \cap L^2(\Omega) $.
\end{itemize}
\end{em}
\end{rem}

\begin{notn}[Generalized weighted total variation; cf.\ {\cite[Section 2]{MS14}}]\label{Note06}
\begin{em}
 \ For any \linebreak $ \varrho \in H^1(\Omega) \cap L^\infty(\Omega) $ and any $ z \in BV(\Omega) \cap L^2(\Omega) $, we define a real-valued Radon measure $ [\varrho |Dz|] \in \mathcal{M}(\Omega) $, as follows:
\begin{equation*}
[\varrho |Dz|](B) := |Dz|_{[\varrho]^+}(B) - |Dz|_{[\varrho]^-}(B) \mbox{ \ for any Borel set $ B \subset \Omega $.}
\end{equation*}
Note that $ [\varrho|\nabla z|](\Omega) $ can be thought of as a generalized version of the total variation of $ z \in BV(\Omega) \cap L^2(\Omega) $ weighted by the possibly sign-changing weight $ \varrho \in H^1(\Omega) \cap L^\infty(\Omega) $. So, hereafter, we simply refer to $ [\varrho|Dz|](\Omega) $ as the \textit{generalized weighted total variation}.
\end{em}
\end{notn}

\begin{rem}\label{Rem.Note06}
\begin{em}
With regard to the generalized weighted total variations, the following facts are verified in \cite[Section 2]{MS14}:
\begin{description}
\item[{(Fact5)}](Strict approximation)
Let $ \varrho \in H^1(\Omega) \cap L^\infty(\Omega) $ and $ z \in BV(\Omega) \cap L^2(\Omega) $ be arbitrary fixed functions, and let $ \{ z_n \, | \, n \in \N \} \subset C^\infty(\overline{\Omega}) $ be any sequence such that
\begin{equation*}
z_n \to z \mbox{ in $ L^2(\Omega) $,  and  strictly in } BV(\Omega), \mbox{ as $ n \to \infty $.}
\end{equation*}
Then
\begin{equation*}
\int_\Omega \varrho |\nabla z_n| \, dx \to \int_\Omega d[\varrho |Dz|] \mbox{ \ as $ n \to \infty $.}
\end{equation*}
\item[{(Fact6)}]For any $ z \in BV(\Omega) \cap L^2(\Omega) $, the mapping
\begin{equation*}
\ds \varrho \in H^1(\Omega) \cap L^\infty(\Omega) \mapsto \int_\Omega d[\varrho |Dz|] \in \R
\end{equation*}
is a linear functional, and moreover, if $ \varphi \in H^1(\Omega) \cap C(\overline{\Omega}) $ and $ \varrho \in H^1(\Omega) \cap L^\infty(\Omega) $, then
\begin{equation*}
\displaystyle \int_\Omega d[\varphi \varrho |Dz|] = \int_\Omega \varphi \, d[\varrho |Dz|].
\end{equation*}
\end{description}
\end{em}
\end{rem}

\begin{notn}[specific classes of functions]\label{Note07}
\begin{em}
Let $ X_0 $ be a Banach space, defined as
\begin{equation*}
X_0 := H^1(\Omega) \times H^1(\Omega) \times BV(\Omega).
\end{equation*}
For arbitrary $ 1 \leq p, q \leq \infty $ and any open interval $ I \subset \R $, we set
\begin{equation*}
L^p(I; BV(\Omega)) := \left\{ \begin{array}{l|l}
\tilde{\theta} & \parbox{8cm}{
$ \tilde{\theta} : I \rightarrow BV(\Omega) $ is measurable for the strict topology of $ BV(\Omega) $, and $ |\tilde{\theta}({}\cdot{})|_{BV(\Omega)} \in L^p(I) $
} 
\end{array} \right\}.
\end{equation*}
As well as, we let:
\begin{equation*}
    L_\mathrm{loc}^p(I; BV(\Omega)) := \left\{ \begin{array}{l|l}
\tilde{\theta} & \parbox{7cm}{
            $ \tilde{\theta} \in L^p(J; BV(\Omega)) $ for an open interval $ J \subset \subset I $, i.e. $ \overline{J} \subset I $
} 
\end{array} \right\}.
\end{equation*}

\noindent
For any $ 1 \leq p \leq \infty $ and any open interval $ I \subset \R $, we note that $ L^p(I; BV(\Omega)) $ and $ L^p(I; X_0) $ are normed spaces, with
\begin{equation*}
|\tilde{\theta}|_{L^p(I; BV(\Omega))} := \bigl| |\tilde{\theta}({}\cdot{})|_{BV(\Omega)} \bigr|_{L^p(I)}, \mbox{ for $ \tilde{\theta} \in L^p(I; BV(\Omega)) $,}
\end{equation*}
and 
\begin{equation*}
\begin{array}{c}
|[\tilde{w}, \tilde{\eta}, \tilde{\theta}]|_{L^p(I; X_0)} := |\tilde{w}|_{L^p(I; H^1(\Omega))} +|\tilde{\eta}|_{L^p(I; H^1(\Omega))} +|\tilde{\theta}|_{L^p(I; BV(\Omega))} 
\\[1ex]
\mbox{for $ [\tilde{w}, \tilde{\eta}, \tilde{\theta}] \in L^p(I; X_0) $,}
\end{array}
\end{equation*}
respectively. 
\end{em}
\end{notn}

Finally, we mention the notion of functional convergence.

\begin{df}[$\Gamma$-convergence; cf.\ \cite{GammaConv}]\label{Def.Gamma}
\begin{em}
Let $ X $ be an abstract Hilbert space, $ \Psi : X \rightarrow (-\infty, \infty] $ be a proper functional, and $ \{ \Psi_n \, | \, n \in \N \} $ be a sequence of proper functionals $ \Psi_n : X \rightarrow (-\infty, \infty] $, $ n \in \N $. We say that $ \Psi_n \to \Psi $ on $ X $, in the sense of $ \Gamma $-convergence \cite{GammaConv}, as $ n \to \infty $ if and only if the following two conditions are fulfilled:
\begin{description}
\item[\textmd{($\gamma$1)}](lower bound) $ \ds \liminf_{n \to \infty} \Psi_n(z_n^\dag) \geq \Psi(z^\dag) $ if $ z^\dag \in X $, $ \{ z_n^\dag \, | \, n \in \N \} \subset X $, and $ z_n^\dag \to z^\dag $ (strongly) in $ X $ as $ n \to \infty $;
\item[\textmd{($\gamma$2)}](optimality) for any $ z^\ddag \in D(\Psi) $, there exists a sequence $ \{ z_n^\ddag  \, | \, n \in \N \} \subset X $ such that $ z_n^\ddag \to z^\ddag $ in $ X $, and $ \Psi_n(z_n^\ddag) \to \Psi(z^\ddag) $, as $ n \to \infty $.
\end{description} 
\end{em}
\end{df}

\section{Assumptions and Main Theorems}

Throughout this paper, we impose the following assumptions. 
\begin{enumerate}
\item[(A0)]
$\Omega \subset \R^{N}$ is a bounded domain of the dimension $N \in \N$, and $\Gamma := \partial\Omega$ is a smooth boundary of $\Omega$. 
Besides, $c \in \R$ is a given constant. 

\item[(A1)]$u \in L^{2}_\mathrm{loc}([0,\infty);L^{2}(\Omega))$ is a given function. 

\item[(A2)]$\alpha_{0} \in W^{1,\infty}_{\rm loc}(\R^{2})$ and $\alpha, \beta \in C^{2}(\R^{2})$ are given functions such that:
\item[~~$\bullet$]$ \alpha $ and $ \beta $ are  convex on $ \R^2 $;  
\item[~~$\bullet$] $ \alpha_\eta(w, 0) \leq 0 $, $ \beta_\eta(w, 0) \leq 0 $, $ \alpha_\eta(w, 1) \geq 0 $, and $ \beta_\eta(w, 1) \geq 0 $, \ for all $ w \in [0, 1] $; 
    
    \item[~~$ \bullet $] $ \ds \inf \bigl[ \alpha_0(\R^{2}) \cup \alpha(\R^{2}) \cup \beta(\R^{2}) \bigr] \geq \delta_{\ast} $, for some $ \delta_* \in (0, 1) $.

\item[(A3)]$ \gamma : \R \to [0, \infty) $ is a proper l.s.c. and convex function such that $ D(\gamma) = [0, 1] $
\item[(A4)] $g \in C^{2}(\R^{2})$ is a function such that
        \begin{equation*}
            \begin{cases}
                g(w, \eta) \geq 0, \mbox{ for all $ [w, \eta] \in \R^2 $,}
                \\[1ex]
                g_\eta(w, 0) \leq 0 \mbox{ and } g_\eta(w, 1) \geq 0, \mbox{ for all $ w \in [0, 1] $.}
            \end{cases}
        \end{equation*}
\item[(A5)] The triplet of initial data $[\v_{0},\theta_{0}] = [w_0, \eta_0, \theta_0] $ belongs to a class $D_{0}$ which is defined by 
\begin{equation*}
D_0 := \left\{ \begin{array}{l|l}
[\tilde{w}, \tilde{\eta}, \tilde{\theta}] \in L^{2}(\Omega)^{3} & 
\parbox{7cm}{
$ \tilde{\theta} \in L^\infty(\Omega) $, and 
$ 0 \leq \tilde{w}, \tilde{\eta} \leq 1 $ a.e. in $ \Omega $}
\end{array} \right\}.
\end{equation*}
\end{enumerate}

For simplicities of descriptions, we prepare the following notations, 
\begin{equation*}
\left\{ 
\begin{array}{l}
G(u; \v) = G(u; w, \eta) := g(w, \eta) + cuw = g(\v) + cuw,
\\[1ex]
[\nabla g](\v) = [\nabla g](w, \eta) := [g_w(w, \eta), g_\eta(w, \eta)] = [g_w(\v), g_\eta(\v)], 
\\[1ex]
\ds [\nabla G](u; \v) = [\nabla G](u; w,\eta) := [g_{w}(w,\eta) + c u, g_{\eta}(w,\eta)] =  [g_{w}(\bm{v}) + c u, g_{\eta}(\bm{v})], 
\end{array}
\right. 
\end{equation*}
and
\begin{equation*}
\begin{array}{c}
\left\{ 
\begin{array}{ll}
\ds [\alpha(\v), \beta(\v)] := [\alpha(w,\eta), \beta(w,\eta)], \vspace{3mm}\\
\ds [ \nabla \alpha ](\v) = [\nabla \alpha](w,\eta) := [\alpha_{w}(w,\eta), \alpha_{\eta}(w,\eta)] = [\alpha_{w}(\v), \alpha_{\eta}(\v)], \vspace{3mm}\\
\ds [ \nabla \beta ](\v) = [\nabla \beta](w,\eta) := [\beta_{w}(w,\eta), \beta_{\eta}(w,\eta)] = [\beta_{w}(\v), \beta_{\eta}(\v)],
\end{array}
\right. 
\end{array}
\end{equation*}
for all $ u \in \R$ and $\v = [w, \eta] \in \R^{2}$. 

Next, for any $\nu \ge 0$ and any $\v = [w,\eta] \in [H^{1}(\Omega) \cap L^{\infty}(\Omega)]^{2}$, we define a convex function $\Phi_{\nu}(\v;{}\cdot\,)$ on $L^{2}(\Omega)$ by 
\begin{equation*}
\theta \in L^{2}(\Omega) \mapsto \Phi_{\nu}(\v;\theta) = \Phi_{\nu}(w,\eta;\theta) := \left\{ \begin{array}{l}
\ds \int_{\Omega} d[\alpha(\v) |D\theta|] + \int_{\Omega} \beta(\v) |D(\nu\theta)|^{2} dx, \vspace{1mm}\\
\ds \hspace{5mm} \mbox{ if } \theta \in BV(\Omega) \mbox{ and } \nu\theta \in H^{1}(\Omega), \vspace{1mm}\\
\ds \infty, \ \ \mbox{ otherwise.}
\end{array}\right. 
\end{equation*}
Note that the above convex function enable us to prescribe the rigorous definition of the free energy $\F_{\nu}$, given on (\ref{free-energy}), as follows, 
\begin{equation}\label{free-energy}
\begin{array}{l}
\ds [\v, \theta] = [w,\eta,\theta] \in L^{2}(\Omega)^{3} \mapsto \F_{\nu}(\v,\theta) = \F_{\nu}(w,\eta,\theta) \vspace{1mm}\\
\ds := \left\{ \begin{array}{l}
\ds \frac{1}{2}\int_{\Omega} |Dw|^{2}dx + \int_{\Omega} \gamma(w) dx 
    + \frac{1}{2} \int_{\Omega} |D\eta|^{2} dx + \int_{\Omega} g(w,\eta) dx + \Phi_{\nu}(w, \eta;\theta) \vspace{1mm}\\
\ds \hspace{5mm} \mbox{ if } \v=[w,\eta] \in [H^{1}(\Omega) \cap L^{\infty}(\Omega)]^{2} \mbox{ and } \theta \in D(\Phi_{\nu}(\v;{}\cdot\,)), \vspace{1mm}\\
\ds \infty,\ \ \mbox{ otherwise. }
\end{array} \right. 
\end{array}
\end{equation}
Now, the solution to the system (S)$_\nu$, for $\nu \ge 0$, is defined as follows. 

\begin{df}\label{df-2}
A triplet $ [\v, \theta] = [w, \eta, \theta] \in L^2_\mathrm{loc}([0, \infty); L^2(\Omega)^3) $ with $ \v = [w, \eta] $ is called  an energy-dissipative solution to $(S)_{\nu}$ or solution to $(S)_{\nu}$ in short, if and only if $ [\v, \theta] $ fulfills the following conditions.
\begin{enumerate}
        \item[(S1)]$ [\bm{v}, \theta] = [w, \eta, \theta] \in C([0, \infty); L^2(\Omega)^3) \cap W_\mathrm{loc}^{1, 2}((0, \infty); L^2(\Omega)^3) 
$
;
                    \\[1ex]
	$\v \in L^{2}_{loc}([0,\infty);H^{1}(\Omega)^{2}) \cap L^{\infty}_{loc}((0,\infty);H^{1}(\Omega)^{2})$; \\[1ex]
	$\theta \in L^{1}_{loc}([0,\infty);BV(\Omega)) \cap L^{\infty}_{loc}((0,\infty); BV(\Omega))$; 
\\[1ex]
                    $ 0 \leq w \leq 1 $, $ 0 \leq \eta \leq 1 $, and $ |\theta| \leq |\theta_0|_{L^\infty(\Omega)} $ a.e. in $ Q $;
                    \\[1ex]
                    $ [\v(0), \theta(0)] = [w(0),\eta(0), \theta(0)]  = [\v_{0}, \theta_{0}]= [w_{0},\eta_{0}, \theta_{0}]  $ in $ L^{2}(\Omega)^{3}$.

\item[(S2)] $ \v = [w, \eta] $ satisfies the following variational forms:
\begin{equation*}
\begin{array}{l}
\ds \left( w_t(t) +g_w(\v(t)) + cu(t), w(t) - \varphi \right)_{L^{2}(\Omega)} + (\nabla w(t), \nabla (w(t)-\varphi))_{L^{2}(\Omega)^{N}} \\[2ex] \qquad
\ds + \int_\Omega d[(w(t)-\varphi)\alpha_{w}(\v(t)) |D \theta(t)|] + \int_{\Omega} (w(t) - \varphi) \beta_{w}(\v(t)) |D(\nu\theta)(t)|^{2} dx \\[2ex] \qquad 
\ds +\int_\Omega \gamma(w(t)) \, dx \le \int_{\Omega}\gamma(\varphi) dx, 
\end{array}
\end{equation*}
for any $ \varphi \in H^1(\Omega) \cap L^{\infty}(\Omega) $ and a.e. $ t \in (0, \infty) $, and
\begin{equation*}
\begin{array}{l}
\ds \left( \eta_t(t) +g_\eta(\v(t)), \psi \right)_{L^{2}(\Omega)} + (\nabla \eta(t), \nabla \psi)_{L^{2}(\Omega)^{N}} \\[1ex] \qquad 
\ds +\int_\Omega d \bigl[ \psi \alpha_{\eta}(\v(t)) |D \theta(t)| \bigr] + \int_{\Omega} \psi \beta_{\eta}(\v(t)) |D(\nu\theta(t))|^{2} dx = 0, 
\end{array}
\end{equation*}
for any $ \psi \in H^1(\Omega) \cap L^{\infty}(\Omega) $ and  a.e. $ t \in  (0, \infty) $.
\item[(S3)] $ \theta $ satisfies the following variational inequality:
\begin{equation*}
\begin{array}{c}
\displaystyle 
(\alpha_0(\v(t)) \theta_t(t), \theta(t) - \omega)_{L^2(\Omega)} + \Phi_{\nu}(\v(t); \theta(t)) \le \Phi_{\nu}(\v(t); \omega),
\end{array}
\end{equation*}
for any $ \omega \in \D(\Phi_{\nu}(\v(t);{}\cdot\,)) $ and  a.e. $ t \in (0, \infty) $. 
\item[(S4)](Energy-dissipation) For any $ {u}^\dag \in L^2(\Omega) $, a function:
    \begin{equation*}
        t \in (0, \infty) \mapsto \F_\nu(\v(t), \theta(t)) -c \bigl( {u}^\dag,  w(t) \bigr)_{L^2(\Omega)} -c^2 \int_0^t |u(t) - {u}^\dag|_{L^2(\Omega)}^2 \, dt,
    \end{equation*}
        coincides with a nonincreasing function belonging to $ L_\mathrm{loc}^1([0, \infty)) \cap BV_\mathrm{loc}((0, \infty)) $. 
\end{enumerate}
\end{df}
\begin{rem}\label{Rem.sols}
\begin{em}
    Two variational forms in (S2) can be reduced to:
\begin{equation}\label{1st.eq}
\begin{array}{l}
\ds (\v_t(t) +[\nabla G](u; \v)(t), \v(t) -\bm{\varpi})_{L^2(\Omega)^2} \vspace{3mm}\\
\ds +(\nabla \v(t), \nabla (\v(t) -\bm{\varpi}))_{L^2(\Omega)^{N \times 2}} 
\ds  +\int_\Omega d \bigl[ (\v(t) -\bm{\varpi}) \cdot [\nabla \alpha](\v(t)) |D \theta(t)| \bigr] \vspace{3mm}\\ 
\ds +\int_\Omega (\v(t) -\bm{\varpi}) \cdot [\nabla \beta](\v(t)) |D(\nu \theta)(t)|^2 \, dx 
\ds +\int_\Omega \gamma(\v(t)) \, dx \leq \int_\Omega \gamma(\bm{\varpi}) \,dx, 
\end{array}
\end{equation}
for any $ \bm{\varpi} = [\varphi, \psi] \in [H^1(\Omega) \cap L^\infty(\Omega)]^2 $ and a.e. $ t \in (0, \infty) $, by using the abbreviation:
\begin{equation*}
	\begin{array}{c}
		\ds \int_\Omega d \bigl[ (\tilde{\bm{{\varpi}}} \cdot \tilde{\v}) |D \tilde{\theta}| \bigr] := \int_\Omega  d \bigl[ \tilde{\varphi} \tilde{w} |D \tilde{\theta}|  \bigr] +\int_\Omega  d \bigl[ \tilde{\psi} \tilde{\eta} |D \tilde{\theta}|  \bigr], 
\end{array}
\end{equation*}
for $ \tilde{\v} = [\tilde{w}, \tilde{\eta}], ~ \tilde{\bm{\varpi}} = [\tilde{\varphi}, \tilde{\psi}] \in [H^1(\Omega) \cap L^\infty(\Omega)]^2 $ and $ \tilde{\theta} \in BV(\Omega) \cap L^2(\Omega) $, and by using the identification 
\begin{equation*}
\gamma(\tilde{\v}) := \gamma (\tilde{w}), 
\end{equation*}
for all $\tilde{\v} = [\tilde{w}, \tilde{\eta}] \in \R^{2}$. 
Furthermore, the variational form in (S3) is equivalent to the following evolution equation:
\begin{equation}\label{3rd.eq}
	\alpha_0(\v(t)) \theta_t(t) +\partial \Phi_{\nu}(\v(t); \theta(t)) \ni 0 \mbox{ in $ L^2(\Omega) $,}
\end{equation}
for a.e. $ t \in (0, \infty)$, governed by the subdifferential $ \partial \Phi_{\nu}(\v(t);{}\cdot\,) \subset L^2(\Omega)^2 $ of the time-dependent convex function $ \Phi_{\nu}(\v(t);{}\cdot\,) $, for $ t \in (0, \infty) $. 
\end{em}
\end{rem}
\medskip 

Now, our Main Theorems are stated as follows. 
\paragraph{Main Theorem 1}{\em 
	Let us assume (A0)--(A5). Then, for any $\nu \ge 0$, the system $(S)_{\nu}$ admits at least one energy-dissipative solution $ [\v, \theta] =  [w, \eta, \theta] \in L^2_\mathrm{loc}([0, \infty); L^2(\Omega)^3) $ with $ \v = [w, \eta] $. 
}

\paragraph{Main Theorem 2}{\em 
In addition to (A0)--(A5), let us assume the following condition. 
\begin{enumerate}
\item[(A6)] There exists a function $u_{\infty} \in L^{2}(\Omega)$ such that $u - u_{\infty} \in L^{2}([0,\infty);L^{2}(\Omega))$.
\end{enumerate}
	Besides, for any $\nu \ge 0$, any solution $ [\v, \theta] =  [w, \eta, \theta]$ to the system (S)$_{\nu}$ with $\v = [w,\eta]$, let us denote by $\omega_{\nu}(\v, \theta)$ the $\omega$-limit set of $[\v,\theta]$ in large time, i.e.: 
\begin{equation*}
\omega_{\nu}(\v, \theta) := \left\{ \begin{array}{l|l}
[w_\infty, \eta_\infty, \theta_\infty] \in L^2(\Omega)^3 & \parbox{7.4cm}{
$ \v_{\infty}=[w_{\infty},\eta_{\infty}] \in [H^{1}(\Omega) \cap L^{\infty}(\Omega)]^{2}, \theta_{\infty} \in BV(\Omega) \cap L^{\infty}(\Omega), \mbox{ and } [\v(t_n), \theta(t_n)] \to [\v_\infty, \theta_\infty] $ in $ L^2(\Omega)^3 $ as $ n \to \infty $, \ for some $ \{ t_n \}_{n = 1}^\infty \subset (0, \infty) $ satisfying $ t_n \uparrow \infty $ as $ n \to \infty $
}
\end{array} \right\}.
\end{equation*}
Then, the following two items hold.
\begin{enumerate}
\item[(O)]$ \omega_{\nu}(\v, \theta) \ne \emptyset $, and $ \omega_{\nu}(\v, \theta) $ is compact in $ L^2(\Omega)^3 $.
\item[(\,I\,)]Any  $ \omega $-limit point, $ [\v_\infty, \theta_\infty] = [w_\infty, \eta_{\infty}, \theta_\infty] \in \omega_{\nu}(\v, \theta) $, fulfills that:
\begin{enumerate}
\item[(i-a)]$ 0 \leq w_{\infty} \leq 1 $, $0 \le \eta_\infty \le 1$, and $ |\theta_\infty| \leq |\theta_0|_{L^\infty(\Omega)} $ a.e. in $ \Omega $;
\item[(i-b)]$ -{\mit \Delta}_{N} \v_\infty + \partial\gamma(\v_{\infty}) +[\nabla G](u_{\infty};\v_\infty) \ni 0 $ in $ L^2(\Omega)^{2} $;
\item[(i-c)]$ \theta_\infty $ is a constant over $ \Omega $, i.e. $ \theta_\infty $ is a global minimizer of the convex function $ \Phi_\nu(\v_\infty;{}\cdot{}) $ on $ L^2(\Omega) $.
\end{enumerate}
\end{enumerate}
}

\section{Approximate problems}
In this section, we introduce approximate problems for the proofs of Main Theorems. The approximation is based on the time discretization method for (\ref{1st.eq})-(\ref{3rd.eq}) with positive constant $\nu$. Therefore, when we consider the approximate problems, we suppose  $ \nu > 0 $, and fix the constant of time-step-size $ h \in (0, 1] $.  

For any $ \nu > 0 $, $\sigma \in (0,1)$, and $\v = [w,\eta] \in [H^{1}(\Omega) \cap L^{\infty}(\Omega)]^{2}$, we define a proper l.s.c. and convex function $ \Phi_\nu^{\sigma}(\v;{}\cdot\,) $ on $ L^2(\Omega) $ by
\begin{equation*}
\ds \theta \in L^2(\Omega) \mapsto \Phi_{\nu}^{\sigma}(\v;\theta) = \Phi_{\nu}^{\sigma}(w,\eta;\theta) := \left\{ 
\begin{array}{ll}
	\multicolumn{2}{l}{
        \ds \int_{\Omega} d[\alpha(\v)|D\theta|_{\sigma}] +\nu^2 \int_\Omega \beta(\v)|D \theta|^2 \, dx,}
	\\[1.5ex]
	& \mbox{if $ \theta \in BV(\Omega) $ and $ \nu  \theta \in H^1(\Omega) $,}\\[2ex]
	\infty, & \mbox{otherwise,}
\end{array}
\right. 
\end{equation*}
and $\Phi_{\nu}^{0}(\v;\theta) := \Phi_{\nu}(\v;\theta)$. 
Here, we take a suitable approximation $\{ |\cdot|_{\sigma}\}_{\sigma \in (0,1)} \subset C^{1}(\R)$ of the Euclidean norm. The precise definition of a suitable approximation is given as follows.
\begin{df}\label{euclregul}
\begin{em}
We say that a collection of functions $\{|\cdot|_\sigma\}_{\sigma\in(0,1)}$ is a suitable approximation to the {Euclidean} norm if the following properties hold.
\begin{description}
    \item [\textmd{(AP1)}] $|\cdot|_\sigma:\R^N\mapsto [0,+\infty)$ is a convex $C^1$ function such that $|0|_\sigma = 0 $ {(and especially, it is differentiable at the origin)}, for all $ \sigma \in (0, 1) $.

\item [\textmd{(AP2)}] There exist bounded functions $ q_{0}\ :\ (0,1) \longrightarrow (0,1]$, $ q_{1} : (0, 1) \longrightarrow [1, \infty) $, $ r_{k}\ :\ (0,1) \longrightarrow [0,\infty)$, $k = 0, 1$, such that:
\begin{equation*}
q_0(\sigma) \to 1, ~ q_1(\sigma) \to 1, ~ r_0(\sigma) \to 0, \ \mbox{ and }\ r_1(\sigma) \to 0, 
\end{equation*}
as $ \sigma \downarrow 0 $, and
\begin{equation*}
\begin{array}{c}
\ds
|\xi|_\sigma \geq q_0(\sigma) |\xi| -r_0(\sigma) \ \mbox{ and }\ |[\nabla |\cdot|_\sigma](\xi)| \leq q_1(\sigma) |\xi|^{r_1(\sigma)},
\end{array}
\end{equation*}
for any $ \xi \in \R^N $ and $ \sigma \in (0, 1) $.
\end{description}
\end{em}
\end{df}
\begin{rem}\label{Rem.Euc}
\begin{em}
Note that (AP1)-(AP2) lead to the following fact:
\begin{center}
$ |\xi|_\sigma \leq [\nabla |\cdot|_\sigma](\xi) \cdot \xi \leq q_1(\sigma)|\xi|^{1 +r_1(\sigma)} $, 
\end{center}
for all $ \xi \in \R^N $ and $ \sigma \in (0, 1) $.
\medskip

    Also, we note that the class of possible regularizations verifying  (AP1)-(AP2) covers a number of standard type regularizations. For instance:
\begin{itemize}
    \item hyperbola type, i.e. $\xi \in \R^N \mapsto \sqrt{|\xi|^2 +\sigma^2} -\sigma$, for $ \sigma \in (0, 1) $;
\item Yosida's regularization, i.e.
    $ \xi \in \R^N \mapsto |\xi|_\sigma := \ds \inf_{\varsigma \in \R^N} \left\{ |\xi| +\frac{\sigma}{2} |\varsigma -\xi|^2 \right\} $, for $ \sigma \in (0, 1) $;
\vspace{-1ex}
\item hyperbolic-tangent type, i.e.
$ \xi \in \R^N \mapsto |\xi|_\sigma := \ds \int_0^{|\xi|} \tanh \frac{\tau}{\sigma} \, d \tau $, for $ \sigma \in (0, 1) $;
\vspace{-1ex}
\item arctangent type, i.e.
$ \xi \in \R^N \mapsto |\xi|_\sigma := \ds \frac{2}{\pi} \int_0^{|\xi|} {\rm Tan}^{-1} \frac{\tau}{\sigma} \, d \tau $, for $ \sigma \in (0, 1) $;
\vspace{-1ex}
\item $ p $-growth type, i.e.
$ \xi \in \R^N \mapsto |\xi|_\sigma := \ds \frac{1}{p(\sigma)} |\xi|^{p(\sigma)} $, for $ \sigma \in (0, 1) $,
with a function $ p : (0, 1) \longrightarrow (1, \infty) $ satisfying $ p(\sigma) \downarrow 1 $ as $  \sigma \downarrow 0 $.
\end{itemize}
\end{em}
\end{rem}
\bigskip

Observe that the convex function $ \Phi_{\nu}^{\sigma}(\v;{}\cdot{}) $ corresponds to a relaxed version of the weighted-total variation $ \Phi_{\nu}(\v;{}\cdot{}) = \Phi_{\nu}^{0}(\v;{}\cdot{}) $.
Additionally, for every $\nu, \sigma \in (0,1)$, we 
define a functional $ \mathscr{F}_{\nu,\sigma}$ on $ L^2(\Omega)^3 $ by letting:
\begin{equation}\label{enrgyF_nu}
\begin{array}{l}
[\v, \theta] = [w, \eta, \theta] \in L^2(\Omega)^3 \mapsto \mathscr{F}_{\nu,\sigma}(\v, \theta) = \mathscr{F}_{\nu,\sigma}(w, \eta, \theta)
\\[1ex]
\qquad := \ds \frac{1}{2} |D \v|_{L^2(\Omega)^{N \times 2}}^2 
    +\int_\Omega \gamma(w) \, dx +\int_\Omega g(w, \eta) \, dx
 + \Phi_\nu^{\sigma}(\v; \theta) 
.
\end{array}
\end{equation}
The above functional $ \mathscr{F}_{\nu,\sigma} $ is a modified version of the free-energy as in \eqref{freeEnrgy}, and the assumptions (A2)--(A4) guarantee the non-negativity of this functional, i.e. $ \mathscr{F}_{\nu,\sigma} \geq 0 $ on $  L^2(\Omega)^3$. Moreover, we note that $\F_{\nu} = \F_{\nu,0}$ for any $\nu \in [0,1)$. 
\bigskip

On this basis, the approximate problem for our system (S)$_\nu$ is denoted by (AP)$_h^{\nu,\sigma}$, and stated as follows. 
\begin{description}
	\item[\textmd{(AP)$_h^{\nu,\sigma}$:}]to find a sequence 
        \begin{equation*}\label{apEq.0} 
            \{ [ \v_i^{\nu,\sigma}, \theta_i^{\nu,\sigma} ] \}_{i = 1}^\infty \subset 
                D_1(\theta_0) := \left\{ 
                    \begin{array}{l|l}
                        [\tilde{\v}, \tilde{\theta}] \in D_0 & \tilde{\theta} \in H^1(\Omega) \mbox{ and } |\tilde{\theta}| \leq |\theta_0|_{L^\infty(\Omega)}
                    \end{array}
                \right\}
        \end{equation*}
with $ \{ \v_i^{\nu,\sigma} \}_{i = 1}^\infty = \{ [w_i^{\nu,\sigma}, \eta_i^{\nu,\sigma}] \}_{i = 1}^\infty  $, which fulfills that
\begin{equation}\label{apEq.2}
\begin{array}{l}
	\ds \frac{1}{h} (\v_i^{\nu,\sigma} -\v_{i -1}^{\nu,\sigma}, \v_i^{\nu,\sigma} -\bm{\varpi})_{L^2(\Omega)^2}+(\nabla \v_i^{\nu,\sigma}, \nabla (\v_i^{\nu,\sigma} -\bm{\varpi}))_{L^2(\Omega)^{N \times 2}}
\\[1ex]
  \ds \qquad +([\nabla G](u_i; \v_i^{\nu,\sigma}), \v_i^{\nu,\sigma} -\bm{\varpi})_{L^2(\Omega)^2}  +\int_\Omega \gamma(\v_i^{\nu,\sigma}) \, dx
	\\[2ex]
	\ds \qquad +\int_\Omega (\v_i^{\nu,\sigma} -\bm{\varpi}) \cdot \bigl( |\nabla \theta_{i -1}^{\nu,\sigma}|_{\sigma} [\nabla \alpha](\v_i^{\nu,\sigma}) +\nu^2 |\nabla \theta_{i -1}^{\nu,\sigma}|^2 [\nabla \beta](\v_i^{\nu,\sigma}) \bigr) \, dx
\\[2ex]
	\ds \qquad  \leq \int_\Omega \gamma(\bm{\varpi}) \, dx, 
\end{array}
\end{equation}
for any $ \bm{\varpi} \in [H^1(\Omega) \cap L^\infty(\Omega)]^2 $,
\begin{equation}\label{apEq.3}
\begin{array}{l}
\ds \frac{1}{h} (\alpha_0(\v_i^{\nu,\sigma}) (\theta_i^{\nu,\sigma} -\theta_{i -1}^{\nu,\sigma}), \omega)_{L^{2}(\Omega)} \vspace{3mm}\\
\ds \hspace{5mm} + (\alpha(\v_{i}^{\nu,\sigma}) [\nabla | \cdot |_{\sigma}](\nabla \theta_{i}^{\nu,\sigma}) + 2\nu^{2} \beta(\v_{i}^{\nu,\sigma})\nabla \theta_{i}^{\nu,\sigma}, \nabla \omega)_{L^{2}(\Omega)} = 0,
\end{array}
\end{equation}
for any $\omega \in H^{1}(\Omega) \cap L^{\infty}(\Omega)$, and any $ i = 1, 2, 3, \dots $, starting from the initial data:

\begin{equation*}
    [\v_{0}^{\nu,\sigma}, \theta_{0}^{\nu,\sigma}] \in D_1(\theta_0)  \mbox{ with } \v_{0}^{\nu,\sigma} = [w_{0}^{\nu,\sigma}, \eta_{0}^{\nu,\sigma}].
\end{equation*}
\end{description}
In the context, for any $ i \in \N $, $ u_i \in L^2(\Omega) $ consists of the components:

\begin{equation*}
	u_i := \frac{1}{h} \int_{(i -1)h}^{ih} [u]_0^{\rm ex}(\tau) \, d\tau \mbox{ in $ L^2(\Omega) $}, 
\end{equation*}
where  $[u]_0^{\rm ex} \in L^{2}(\R; L^{2}(\Omega))$ is the zero-extensions of $u$. 
\medskip

Now, before the proof of Main Theorem 1, it will be needed to verify the following lemmas.
\begin{lem}[Solvability of the approximate problem]\label{thm1}
    There exists a small constant $ h_1^\dagger \in (0, 1] $ such that if $\nu,\sigma > 0$ and $ h \in (0,  h_1^\dagger] $, then the approximate problem (AP)$^{\nu,\sigma}_{h}$ admits a unique solution $ \{ [\v_{i}^{\nu,\sigma}, \theta_{i}^{\nu,\sigma}] \}_{i = 1}^\infty \subset D_1(\theta_0) $, and moreover,
\begin{equation}\label{AP_diss}
\begin{array}{c}
    \ds \frac{1}{2h}|\v_{i}^{\nu,\sigma} - \v_{i-1}^{\nu,\sigma}|_{L^{2}(\Omega)^{2}}^{2}  + \frac{1}{h} \bigl| {\ts \sqrt{\alpha_{0}(\v_{i}^{\nu,\sigma})}(\theta_{i}^{\nu,\sigma}-\theta_{i-1}^{\nu,\sigma})} \bigr|_{L^{2}(\Omega)}^{2} + \F_{\nu,\sigma}(\v_{i}^{\nu,\sigma},\theta_{i}^{\nu,\sigma}) \vspace{3mm}\\
\ds + c({u}^\dag,w_{i}^{\nu,\sigma})_{L^{2}(\Omega)}
\ds 
\ds \le \F_{\nu,\sigma}(\v_{i-1}^{\nu,\sigma},\theta_{i-1}^{\nu,\sigma}) + c({u}^\dag,w_{i-1}^{\nu,\sigma})_{L^{2}(\Omega)} + c^{2}h|u_{i}-{u}^\dag|_{L^{2}(\Omega)}^{2}, 
\end{array}
\end{equation}
for $ i = 1, 2, 3, \dots$ and any ${u}^\dag \in L^{2}(\Omega)$, and 
\begin{equation}\label{AP_diss2}
\begin{array}{c}
    \ds \frac{1}{2} \sum_{i=1}^{m} i | \v_{i}^{\nu,\sigma} - \v_{i-1}^{\nu,\sigma}|_{L^{2}(\Omega)^{2}}^{2} + \sum_{i=1}^{m} i \bigl| {\ts \sqrt{\alpha_{0}(\v_{i}^{\nu,\sigma})}(\theta_{i}^{\nu,\sigma} - \theta_{i-1}^{\nu,\sigma})} \bigr|_{L^{2}(\Omega)}^{2} \vspace{3mm}\\
\ds + mh \F_{\nu,\sigma}(\v_{m}^{\nu,\sigma},\theta_{m}^{\nu,\sigma}) + cmh ({u}^\dag, w_{m}^{\nu,\sigma})_{L^{2}(\Omega)} \vspace{3mm}\\
\ds \le h \sum_{i=1}^{m} \F_{\nu,\sigma}(\v_{i-1}^{\nu,\sigma},\theta_{i-1}^{\nu,\sigma}) + ch \sum_{i=1}^{m} ({u}^\dag,w_{i-1}^{\nu,\sigma})_{L^{2}(\Omega)} + c^{2}h^{2}\sum_{i=1}^{m}i|u_{i}-{u}^\dag|_{L^{2}(\Omega)}^{2}, 
\end{array}
\end{equation}
for any $m \in \N$ and any ${u}^\dag \in L^{2}(\Omega)$. 
\end{lem}
\bigskip

\noindent
\textbf{Proof.} 
By way of a slight modification of the proof of \cite[Theorem 1]{SWY14}, the existence and uniqueness of approximate solutions are verified. 

To show the inequality (\ref{AP_diss}), we set ${\bm \varpi} = \vv_{i-1}$ in (\ref{apEq.2}). By using (A2) and Young's inequality, we have 
\begin{align*}
    \ds \frac{1}{h} | &{} \vv_{i} - \vv_{i-1}|_{L^{2}(\Omega)^{2}}^{2} + \frac{1}{2}|\nabla \vv_{i}|_{L^{2}(\Omega)^{2N}}^{2}  + \int_{\Omega} \gamma(\vv_{i}) dx 
    \\
    & \ds + \int_{\Omega} \alpha(\vv_{i}) |\nabla \thetaa_{i-1}|_{\sigma} dx + \nu^{2} \int_{\Omega} \beta(\vv_{i}) |\nabla \thetaa_{i-1}|^{2} dx 
    \\
    & + \int_{\Omega} [\nabla g](\vv_{i}) \cdot (\vv_{i} - \vv_{i-1}) dx \,+c \int_{\Omega} u_{i}(\ww_{i} - \ww_{i-1}) dx
    \\
    \le &~ \ds \frac{1}{2}|\nabla \vv_{i-1}|_{L^{2}(\Omega)^{2N}}^{2} + \int_{\Omega} \gamma(\vv_{i-1}) dx 
    \\
    & + \int_{\Omega} \alpha(\vv_{i-1}) |\nabla \thetaa_{i-1}|_{\sigma} dx + \nu^{2} \int_{\Omega} \beta(\vv_{i-1}) |\nabla \thetaa_{i-1}|^{2} dx.
\end{align*}
According to (A4), we note that 
\begin{equation*}
\begin{array}{l}
\ds g(\vv_{i-1}) \ge g(\vv_{i}) + [\nabla g](\vv_{i})\cdot(\vv_{i-1}-\vv_{i}) - \frac{1}{2}|g|_{C^{2}([0,1]^{2})}|\vv_{i-1}-\vv_{i}|^{2}, 
\end{array}
\end{equation*}
and hence 
\begin{equation}\label{kenHelp01}
    \begin{array}{l}
        \ds \int_{\Omega} [\nabla g](\vv_{i}) \cdot (\vv_{i} - \vv_{i-1}) dx \vspace{3mm}\\
        \ds \ge \int_{\Omega} g(\vv_{i})dx - \int_{\Omega} g(\vv_{i-1})dx - \frac{1}{2} |g|_{C^{2}([0,1]^{2})} |\vv_{i} - \vv_{i-1}|_{L^{2}(\Omega)^{2}}^{2}.
    \end{array}
\end{equation}
Also, using Young's inequality again, for any ${u}^\dag \in L^{2}(\Omega)$, it holds that:
\begin{equation*}
\begin{array}{l}
\ds c\int_{\Omega} u_{i}(\ww_{i} - \ww_{i-1}) dx 
\ds = c\int_{\Omega} {u}^\dag(\ww_{i} - \ww_{i-1}) dx + c\int_{\Omega} (u_{i}-{u}^\dag)(\ww_{i} - \ww_{i-1}) dx \vspace{3mm}\\
\ds \ge c\int_{\Omega} {u}^\dag(\ww_{i} - \ww_{i-1}) dx - c^{2}h\int_{\Omega} |u_{i}-{u}^\dag|^{2} dx - \frac{1}{4h} \int_{\Omega} |\ww_{i} - \ww_{i-1}|^{2} dx.
\end{array}
\end{equation*}
Hence, we can get the following inequality
\begin{align}\label{inq4-1}
    & \ds \left( \frac{3}{4h}-\frac{|g|_{C^{2}([0,1]^{2})}}{2} \right) |\vv_{i} - \vv_{i-1}|_{L^{2}(\Omega)^{2}}^{2} + \frac{1}{2}|\nabla \vv_{i}|_{L^{2}(\Omega)^{2N}}^{2} \vspace{2mm}
    \nonumber
    \\
    & \qquad \ds + \int_{\Omega} \gamma(\vv_{i}) dx + \int_{\Omega} g(\vv_{i})dx 
    +\Phi^{\nu, \sigma}(\v_i^{\nu, \sigma}; \theta_{i -1}^{\nu, \sigma})
    + c\int_{\Omega} {u}^\dag \ww_{i} dx
    \nonumber
    \\
    & \ds \le \frac{1}{2}|\nabla \vv_{i-1}|_{L^{2}(\Omega)^{2N}}^{2} + \int_{\Omega} \gamma(\vv_{i-1}) dx + \int_{\Omega} g(\vv_{i-1})dx
    \\
    & \qquad \ds 
    +\Phi^{\nu, \sigma}(\vv_{i -1}; \theta_{i -1}^{\nu, \sigma}) + c\int_{\Omega} {u}^\dag \ww_{i-1}dx + c^{2}h\int_{\Omega} |u_{i}-{u}^\dag|^{2} dx. 
    \nonumber
\end{align}
On the other hand, we set $\omega = \thetaa_{i} - \thetaa_{i-1}$ in (\ref{apEq.3}). Then, Remark \ref{Rem.Euc} and Young's inequality yield 
\begin{equation}\label{inq4-2}
    \ds \frac{1}{h} \bigl| {\ts \sqrt{\alpha_{0}(\vv_{i})}(\thetaa_{i} - \thetaa_{i-1})} \bigr|_{L^{2}(\Omega)}^{2} + \Phi^{\nu,\sigma}(\vv_{i}; \thetaa_{i}) \le \Phi^{\nu,\sigma}(\vv_{i}; \thetaa_{i-1}). 
\end{equation}
Here, we set
\begin{equation}\label{h1}
h^{\dagger}_{1} := \frac{1}{2(1 \vee |g|_{C^{2}([0,1]^{2})})}.
\end{equation}
Since 
\begin{equation*}\label{h}
    \frac{3}{4h} - \frac{|g|_{C^{2}([0,1]^{2})}}{2} > \frac{1}{2h}, 
\end{equation*}
for $0 < h < h^{\dagger}_{1}$, the desired inequality (\ref{AP_diss}) is obtained by taking the sum of (\ref{inq4-1}) and (\ref{inq4-2}). 

To prove (\ref{AP_diss2}), we multiply both sides of (\ref{AP_diss}) by $ih$, and take a summation of the inequality from $1$ to $m \in \N$. 
\hfill $\Box$ 
\bigskip

\begin{lem}\label{lem4-1}
Let $\nu, \sigma \in (0,1)$, let $h_{1}^{\dagger}$ be the constant in (\ref{h1}), let $h \in (0,h_{1}^{\dagger})$ be an arbitrary time-step, and let $\{\vv_{i}, \thetaa_{i}\}$ be the solution to (AP)$_{h}^{\nu,\sigma}$ with initial data $[\v_{0}^{\nu,\sigma},\theta_{0}^{\nu,\sigma}] \in D_{1}(\theta_{0})$. Under assumptions (A0)-(A5), there exist $\nu_{\ast} \in (0,1)$ and positive constants $A_{\ast}$, $B_{\ast}$, $C_{\ast}$, depending only on $\Omega$, $\alpha_{0}$, $\alpha$, $\beta$, $g$, $\gamma$, and $\theta_{0}$, such that if $h \in (0,h_{1}^{\dagger})$ and $\nu \in (0,\nu_{\ast})$, then the approximate solution $\{\vv_{i}, \thetaa_{i}\}$ satisfies the following energy inequality: 
\begin{equation}\label{lem1-inq}
\begin{array}{l}
\ds \frac{1}{2}(|\vv_{m} -\w_{0}|_{L^{2}(\Omega)^{2}}^{2} + A_{\ast} |\thetaa_{m} -\omega_{0}|_{L^{2}(\Omega)}^{2}) + \frac{B_* h}{2} \sum_{i = 1}^{m} \mathscr{F}_{\nu, \sigma}(\vv_{i-1}, \thetaa_{i-1})\vspace{3mm}\\
\ds \le \ds \frac{1}{2}(|\v_{0}^{\nu,\sigma}-\w_{0}|_{L^{2}(\Omega)^{2}}^{2} + A_{\ast} |\theta_{0}^{\nu,\sigma}-\omega_{0}|_{L^{2}(\Omega)}^{2}) +\frac{h}{B_*} \mathscr{F}_{\nu, \sigma}(\v_{0}^{\nu}, \theta_{0}^{\nu})\vspace{3mm}\\
\ds\hspace{5mm} +m h C_*(1 +|\w_0|_{H^1(\Omega)^{2}}^2 +|\omega_0|_{H^1(\Omega)}^2) + \frac{c^{2}h}{2} \sum_{i=1}^{m} |u_{i}|_{L^{2}(\Omega)}^{2},
\end{array}
\end{equation}
for any $ m \in \N $ and any $ [\w_0, \omega_0]\in D_1(\theta_0)$.
\end{lem}
\textbf{Proof.}
Let us fix a pair (triplet) of functions:
    \begin{equation*}\label{ken00}
        [\w_0, \omega_0] \in D_1(\theta_0) \mbox{ with } \w_0 = [\tilde{w}_0, \tilde{\eta}_0] \in H^1(\Omega)^2,
    \end{equation*}
    and fix a time-step $i \in \N$. Also, we define a large constant $R_{\ast} > 0$ by 
\begin{equation*}
\begin{array}{l}
\ds R_* := \bigl[ (1 +|\alpha_0|_{W^{1, \infty}((0, 1)^{2})})(1 +|\alpha|_{C^1([0, 1]^{2})})(1 + |\beta|_{C([0,1]^{2})}) 
\vspace{3mm}\\
    \hspace{12mm} \ds \cdot (1 + |\gamma|_{L^{\infty}(0, 1)})(1 +|g|_{W^{2, \infty}((0, 1)^{2})})(1 +|\theta_0|_{L^\infty(\Omega)})(1 +\mathscr{L}^N(\Omega)) \bigr]^{2} / \delta_{\ast}^{4}.
\end{array}
\end{equation*}

First let us set ${\bm \varpi} = {\bm w_{0}} = [\widetilde{w}_{0}, \widetilde{\eta}_{0}]$ in (\ref{apEq.2}). Then, using \eqref{kenHelp01}, (A2), and Young's inequality, we have 
\begin{align*}
    & \ds \frac{1}{2h}(|\vv_{i} - \w_{0}|_{L^{2}(\Omega)^{2}}^{2} - |\vv_{i-1} - \w_{0}|_{L^{2}(\Omega)^{2}}^{2}) + \frac{1}{2}|\nabla \vv_{i}|_{L^{2}(\Omega)^{2N}}^{2} 
    \\
    & \ds + \int_{\Omega} \gamma(w_{i}^{\nu, \sigma})dx + 
    \int_{\Omega} g(\vv_{i})dx 
    - 3|g|_{W^{2,\infty}((0,1)^{2})} \mathscr{L}^{N}(\Omega) 
    \\
    & \ds + \int_{\Omega} |\nabla \thetaa_{i-1}|_{\sigma} (\alpha(\vv_{i}) - \alpha(\w_{0}))dx + \nu^{2} \int_{\Omega} |\nabla \thetaa_{i-1}|^{2} (\beta(\vv_{i}) - \beta(\w_{0}))dx
    \\
    & \ds  +c(u_{i}, \ww_{i} - \widetilde{w}_{0})_{L^{2}(\Omega)} \le \frac{1}{2}|\nabla \w_{0}|_{L^{2}(\Omega)^{2N}}^{2} + \int_{\Omega} \gamma(\tilde{w}_{0}) dx.
\end{align*}
By (A2) and H\"{o}lder's inequality, it is deduced that 
\begin{equation}\label{inq4-4}
\begin{array}{l}
\ds \frac{1}{2h}(|\vv_{i} - \w_{0}|_{L^{2}(\Omega)^{2}}^{2} - |\vv_{i-1} - \w_{0}|_{L^{2}(\Omega)^{2}}^{2}) + \frac{1}{2}|\nabla \vv_{i}|_{L^{2}(\Omega)^{2N}}^{2} \vspace{3mm}\\
\ds + \int_{\Omega} \gamma(w_{i}^{\nu,\sigma}) dx + 
\int_{\Omega} g(\vv_{i})dx
\vspace{3mm}\\
\ds + \frac{\delta_{\ast}}{|\alpha|_{C([0,1]^{2})}} \int_{\Omega} \alpha(\vv_{i-1})|\nabla \thetaa_{i-1}|_{\sigma} dx - |\alpha|_{C([0,1]^{2})} \int_{\Omega} |\nabla \thetaa_{i-1}|_{\sigma} dx\vspace{3mm}\\
\ds   + \nu^{2} \frac{\delta_{\ast}}{|\beta|_{C([0,1]^{2})}}\int_{\Omega} \beta(\vv_{i-1})|\nabla \thetaa_{i-1}|^{2} dx - \nu^{2} |\beta|_{C([0,1]^{2})} \int_{\Omega} |\nabla \thetaa_{i-1}|^{2} dx \vspace{3mm}\\
\ds \le 4R_*(1 + |\w_{0}|_{H^{1}(\Omega)^{2}}^{2}) + \frac{c^{2}}{2} |u_{i}|_{L^{2}(\Omega)}^{2}.
\end{array}
\end{equation}

Next, we take $\omega := (\theta_{i}^{\nu, \sigma} - \omega_{0}) / \alpha_{0}(\v_{i}^{\nu, \sigma})$ as the test function in (\ref{apEq.3}). Then, from (A2) and Definition \ref{euclregul}, we see that 
\begin{align}\label{inq4-3}
    \ds \frac{1}{2h} & (|\thetaa_{i} - \omega_{0}|_{L^{2}(\Omega)}^{2} - |\thetaa_{i-1} - \omega_{0}|_{L^{2}(\Omega)}^{2}) 
    +\int_\Omega \frac{\alpha(\vv_i)}{\alpha_0(\vv_i)} [\nabla |\cdot|_\sigma](\nabla\thetaa_{i}) \cdot \nabla (\thetaa_i -\omega_0) \, dx
    \nonumber
    \\
    & + 2 \nu^{2} \int_{\Omega} \frac{\beta(\vv_{i})}{\alpha_{0}(\vv_{i})} \nabla \thetaa_{i} \cdot \nabla (\thetaa_{i} - \omega_{0}) dx
    \nonumber
    \\
    \le & \ds \int_{\Omega} \frac{\alpha(\vv_{i}) (\thetaa_{i} - \omega_{0}) }{\alpha_{0}(\vv_{i})^{2}} [\nabla|\cdot|_{\sigma}](\nabla\thetaa_{i}) \cdot \nabla \alpha_{0}(\vv_{i}) dx
    \\
    & \ds \hspace{5mm} + 2 \nu^{2} \int_{\Omega} \frac{\beta(\vv_{i})(\thetaa_{i} - \omega_{0})}{\alpha_{0}(\vv_{i})^{2}} \nabla \thetaa_{i} \cdot \nabla \alpha_{0}(\vv_{i}) dx,  
    \nonumber
\end{align}
and by using Remark \ref{Rem.Euc} and Young's inequality, we also have 
\begin{subequations}\label{inq4-3_ab}
    \begin{align}\label{inq4-3_c}
    \ds \int_{\Omega} & \frac{\alpha(\vv_{i})}{\alpha_{0}(\vv_{i})} [\nabla|\cdot|_{\sigma}](\nabla\thetaa_{i}) \cdot \nabla (\thetaa_{i} - \omega_{0}) dx \vspace{3mm}
    \nonumber
    \\
    & \ds \ge \int_{\Omega} \frac{\alpha(\vv_{i})}{\alpha_{0}(\vv_{i})} |\nabla\thetaa_{i}|_{\sigma} dx - \int_{\Omega} \frac{\alpha(\vv_{i})}{\alpha_{0}(\vv_{i})} q_{1}(\sigma)|\nabla\omega_{0}|^{1 + r_{1}(\sigma)} dx
    \\
        & \ds \ge \frac{\delta_*}{|\alpha_0|_{C([0, 1]^2)}} \int_\Omega |\nabla \thetaa_i|_\sigma \, dx -\frac{|\alpha|_{C([0, 1]^2)}}{\delta_*} \int_\Omega q_{1}(\sigma) |\nabla\omega_{0}|^{1 + r_{1}(\sigma)} dx,
    \nonumber
\end{align}
  
\begin{align}
    \ds 2 \nu^2 \int_{\Omega} &~ \frac{\beta(\vv_{i})}{\alpha_{0}(\vv_{i})} \nabla \thetaa_{i} \cdot \nabla (\thetaa_i -\omega_{0}) dx  
    \nonumber
    \\
    \ge &~ 2 \nu^2 \cdot \frac{1}{|\alpha_0|_{C([0, 1]^2)}} \int_\Omega \beta(\vv_i) |\nabla \thetaa_i|^2 dx 
    \nonumber
    \\
    & \quad -2 \nu^2 \int_\Omega \left( {\ts \sqrt{\beta(\vv_i)}|\nabla \thetaa_i|} \right) \left( \rule{-1pt}{14pt} \right. \frac{1}{\delta_*} {\ts \sqrt{\beta(\vv_i)} |\nabla \omega_0|} \left. \rule{-2pt}{14pt} \right) dx
    \nonumber
    \\
    \geq & \frac{3\nu^2}{2 |\alpha_0|_{C([0, 1]^2)}} \bigl| {\ts \sqrt{\beta(\vv_i)} \nabla \thetaa_i} \bigr|_{L^2(\Omega)^N}^2 
    \nonumber
    \\
    & \quad -\frac{2 \nu^2 |\alpha_0|_{C([0, 1]^2)}}{\delta_*^2}\bigl| {\ts \sqrt{\beta(\vv_i)} \nabla \omega_0} \bigr|_{L^2(\Omega)^N}^2,
    \label{inq4-3_a}
\end{align}
and 
\begin{align}\label{inq4-3_b}
    & \ds 2 \nu^2 \int_{\Omega} \frac{\beta(\vv_{i})(\thetaa_{i} - \omega_{0})}{\alpha_{0}(\vv_{i})^{2}} \nabla \thetaa_{i} \cdot \nabla \alpha_{0}(\vv_{i}) dx 
    \nonumber
    \\
    \ds \le & 2 \nu^2 |\alpha_{0}|_{C([0,1]^{2})} \cdot \frac{|\thetaa_i -\omega_0|_{L^\infty(\Omega)}^2 |\beta|_{C([0,1]^{2})}}{\delta_{\ast}^{4}} \cdot \bigl| \, |[ \nabla \alpha_{0}] | \, \bigr|_{L^{\infty}((0,1)^{2})}^{2} |\nabla \vv_{i}|_{L^{2}(\Omega)^{2N}}^{2}
    \nonumber
    \\
    & \ds \hspace{5mm} + \frac{\nu^2}{2|\alpha_{0}|_{C([0,1]^{2})}} \bigl| {\ts \sqrt{\beta(\vv_{i})}\nabla\thetaa_{i}} \bigr|_{L^{2}(\Omega)^{N}}^{2}
    \nonumber
    \\
    \ds \le & 16 \nu^2 |\alpha_{0}|_{C([0,1]^{2})} \cdot \frac{|\theta_0|_{L^\infty(\Omega)}^2|\alpha_{0}|_{W^{1, \infty}((0,1)^{2})}^{2} |\beta|_{C([0,1]^{2})}}{\delta_{\ast}^{4}} |\nabla \vv_{i}|_{L^{2}(\Omega)^{2N}}^{2}
    \\
    & \ds \hspace{5mm} + \frac{\nu^2}{2|\alpha_{0}|_{C([0,1]^{2})}} \ts\bigl| \sqrt{\beta(\vv_{i})}\nabla\thetaa_{i} \bigr|_{L^{2}(\Omega)^{N}}^{2}.
    \nonumber
\end{align}
\end{subequations}
Here, let us set:
\begin{equation}\label{A-ast}
A_{\ast} := \frac{2|\alpha_{0}|_{C([0,1]^{2})} \max\{|\alpha|_{C([0,1]^{2})}, |\beta|_{C([0,1]^{2})}\}}{\delta_{\ast
}} ~ (\le 2\delta_{\ast
}R_{\ast}^{1/2}). 
\end{equation}
Then, by Remark \ref{Rem.Euc} and Young's inequality, we compute that
\begin{align}\label{inq4-3_d}
    \ds \int_{\Omega} & \frac{\alpha(\vv_{i}) (\thetaa_{i} - \omega_{0}) }{\alpha_{0}(\vv_{i})^{2}} [\nabla|\cdot|_{\sigma}](\nabla\thetaa_{i}) \cdot \nabla \alpha_{0}(\vv_{i}) dx
    \nonumber
    \\
    & \ds \le 8A_{\ast}R_{\ast}\int_{\Omega} q_{1}(\sigma)^{2}|\nabla\thetaa_{i}|^{2r_{1}(\sigma)} dx  + \frac{1}{8A_{\ast}} |\nabla \vv_{i}|_{L^{2}(\Omega)^{2}}^{2}
\end{align}
Moreover, multiplying the both sides of \eqref{inq4-3} by $ A_* $, it is inferred from \eqref{inq4-3_ab}--\eqref{inq4-3_d} that:
\begin{align}\label{inq4-4.5}
    \ds \frac{A_{\ast}}{2h} & (|\thetaa_{i} - \omega_{0}|_{L^{2}(\Omega)}^{2} - |\thetaa_{i-1} - \omega_{0}|_{L^{2}(\Omega)}^{2}) + 2 \mbox{max}\{|\alpha|_{C([0,1]^{2})}, |\beta|_{C([0,1]^{2})}\} \bigl| |\nabla\thetaa_{i}|_{\sigma} \bigr|_{L^{1}(\Omega; \R^N)} \nonumber
    \\
    & \ds + \frac{2 \nu^{2}\max\{|\alpha|_{C([0,1]^{2})},|\beta|_{C([0,1]^{2})}\}}{\delta_{\ast}} \bigl| {\ts \sqrt{\beta(\vv_{i})} \nabla \thetaa_{i}} \bigr|_{L^{2}(\Omega)^{N}}^{2} \nonumber
    \\
    & \ds - \left( \frac{1}{8} +16\nu^{2}|\alpha_{0}|_{C([0,1]^{2})}A_{\ast}R_{\ast} \right) |\nabla \vv_{i}|_{L^{2}(\Omega)^{2N}}^{2} \nonumber
    \\
    \ds \le & 2\nu^{2} A_{\ast}\frac{|\alpha_{0}|_{C([0,1]^{2})}}{\delta_{\ast}^{2}} \bigl| {\ts \sqrt{\beta(\vv_{i})} \nabla \omega_{0}} \bigr|_{L^{2}(\Omega)^{N}}^{2} \nonumber
    \\
    & \ds + \frac{A_{\ast} |\alpha|_{C([0,1]^{2})}}{\delta_{\ast}} \int_{\Omega} q_{1}(\sigma)|\nabla\omega_{0}|^{1 + r_{1}(\sigma)} dx + 8A_{\ast}^{2}R_{\ast}\int_{\Omega} q_{1}(\sigma)^{2}|\nabla\thetaa_{i}|^{2r_{1}(\sigma)} dx. 
\end{align}

Now, we set a constant $ \nu_* \in (0, 1) $ so small to satisfy that:
\begin{subequations}\label{ken10}
\begin{equation}\label{nu-ast}
\ds 0 < \nu_{\ast}^{2} < \min \left\{ \frac{1}{128|\alpha_{0}|_{C([0,1]^{2})}A_{\ast}R_{\ast}},\ \frac{1}{2} \right\}, 
\end{equation}
and
\begin{equation}
\left \{ 
\begin{array}{l}\label{q-r}
\ds \frac{3}{4} \le q_{0}(\sigma) \le 1 \le q_{1}(\sigma) \le \frac{5}{4},
\\[1ex]
\ds 0 \le r_{0}(\sigma) \le \frac{1}{4},\ \ 0 \le r_{1}(\sigma) \le \frac{1}{4},
\end{array} \mbox{for all $ \sigma \in (0, \nu_*) $.}
\right. 
\end{equation}
\end{subequations}
Here, we use the following type Young's inequality: for arbitrary $a, b \geq 0 $,  $\widehat{\varepsilon} \in (0, 1)$, and $ 1 \leq q \leq 2 \leq p < \infty$ with $\frac{1}{p} + \frac{1}{q} =1$,
\begin{equation}\label{young}
    ab = \bigl( (p \widehat{\varepsilon})^{\frac{1}{p}} a \bigr) \bigl( (p \widehat{\varepsilon})^{-\frac{1}{p}} b \bigr) \leq \widehat{\varepsilon} a^p +\frac{1}{q} ( p \widehat{\varepsilon})^{-\frac{q}{p}} b^q \le \widehat{\varepsilon} a^{p} + \frac{1}{\widehat{\varepsilon}}b^{q};
\end{equation}
to estimate the third term on the right-hand side of \eqref{inq4-4.5}. Indeed, by letting:
\begin{equation*}
    a = |\nabla \thetaa_i|^{2 r_1(\sigma)}, ~ b = 1, \mbox{ and } p = \frac{1}{2 r_1(\sigma)} \geq 2,
\end{equation*}
to apply \eqref{young}, the third term on the right-hand side of \eqref{inq4-4.5} is estimated as follows:
\begin{equation}\label{inq4-5}
\begin{array}{l}
\ds 8A_{\ast}^{2}R_{\ast}\int_{\Omega} q_{1}(\sigma)^{2}|\nabla\thetaa_{i}|^{2r_{1}(\sigma)} dx 
\ds \le \frac{25}{2} A_{\ast}^{2}R_{\ast} \int_{\Omega} \left( \widehat{\varepsilon} |\nabla \thetaa_{i}| + \frac{1}{\widehat{\varepsilon}} \right) dx  \vspace{3mm}\\
\ds \le \frac{50R_{\ast}^{2}\delta_{\ast
    }^{2}\widehat{\varepsilon}}{q_{0}(\sigma)} \bigl( \bigl| |\nabla\thetaa_{i}|_{\sigma} \bigr|_{L^{1}(\Omega; \R^N)} + r_{0}(\sigma) \mathscr{L}^{N}(\Omega) \bigr)+ \frac{50R_{\ast}^{2}\delta_{\ast
}^{2}
}{\widehat{\varepsilon}}\mathscr{L}^{N}(\Omega),
\end{array}
\end{equation}
by (\ref{q-r}), 
(AP2), and (\ref{A-ast}). Then, for any $0 < \widehat{\varepsilon} < 1$, it follows that 
\begin{align*}
    \ds \frac{A_{\ast}}{2h} ( & |\thetaa_{i} - \omega_{0}|_{L^{2}(\Omega)}^{2} - |\thetaa_{i-1} - \omega_{0}|_{L^{2}(\Omega)}^{2}) - \frac{1}{4} |\nabla \v_{i}^{\nu}|_{L^{2}(\Omega)^{2N}}^{2}
    \\
    & + 2\mbox{max}\{|\alpha|_{C([0,1]^{2})}, |\beta|_{C([0,1]^{2})}\} \bigl| |\nabla\thetaa_{i}|_{\sigma} \bigr|_{L^{1}(\Omega; \R^N)}
    \\
    & \ds + \frac{2\nu^{2}\max\{|\alpha|_{C([0,1]^{2})},|\beta|_{C([0,1]^{2})}\}}{\delta_{\ast
    }} \bigl| {\ts \sqrt{\beta(\vv_{i})} \nabla \thetaa_{i}} \bigr|_{L^{2}(\Omega)^{N}}^{2}
    \\[2ex]
    \le &~ \ds 4\nu^{2}R_{\ast}|\nabla \omega_{0}|_{L^{2}(\Omega)^{N}}^{2} 
    \\
    & + \frac{5A_{\ast} |\alpha|_{C([0,1]^{2})}}{4\delta_{\ast
}} \left( \frac{1+r_{1}(\sigma)}{2}|\nabla\omega_{0}|_{L^{2}(\Omega)^{N}}^{2} + \frac{1-r_{1}(\sigma)}{2} \mathscr{L}^{N}(\Omega) \right)
    \\
    & \ds + \frac{50R_{\ast}^{2}\delta_{\ast
    }^{2}\widehat{\varepsilon}}{q_{0}(\sigma)} \bigl( \bigl||\nabla\thetaa_{i}|_{\sigma}\bigr|_{L^{1}(\Omega; \R^N)} + r_{0}(\sigma) \mathscr{L}^{N}(\Omega)\bigr)+ \frac{50R_{\ast}^{2}\delta_{\ast
}^{2}
}{\widehat{\varepsilon}}\mathscr{L}^{N}(\Omega)
    \\[2ex]
    \le &~ \ds \frac{200R_{\ast}^{2}\delta_{\ast}^{2}\widehat{\varepsilon}}{3} \bigl||\nabla\thetaa_{i}|_{\sigma} \bigr|_{L^{1}(\Omega; \R^N)} + 50 R_{\ast}^{3} \left( \delta_{\ast
}^{2}\widehat{\varepsilon} + \frac{\delta_{\ast
}^{2}}{\widehat{\varepsilon}} + 1 \right) (1 + |\nabla \omega_{0}|_{L^{2}(\Omega)^{N}}^{2})
    \\[1ex]
    \le & \ds 10^{2} R_{\ast}^{2} \widehat{\varepsilon} \bigl||\nabla\thetaa_{i}|_{\sigma}\bigr|_{L^{1}(\Omega; \R^N)} + 10^{2} R_{\ast}^{3} \left( \widehat{\varepsilon} + \frac{\delta_{\ast
}^{2}}{\widehat{\varepsilon}} + 1 \right) (1 + |\nabla \omega_{0}|_{L^{2}(\Omega)^{N}}^{2}), 
\end{align*}
by (A2), 
(\ref{A-ast})-(\ref{q-r}), (\ref{inq4-5}), and Young's inequality. If we choose 
\begin{equation*}
\widehat{\varepsilon} := \frac{1}{10^{2}R_{\ast}^{2}}|\alpha|_{C([0,1]^{2})}, 
\end{equation*}
we can see that 
\begin{align}\label{inq4-6}
    \ds \frac{A_{\ast}}{2h} & (|\thetaa_{i} - \omega_{0}|_{L^{2}(\Omega)}^{2} - |\thetaa_{i-1} - \omega_{0}|_{L^{2}(\Omega)}^{2}) - \frac{1}{4} |\nabla \vv_{i}|_{L^{2}(\Omega)^{2N}}^{2} 
    \nonumber
    \\
    &~ + \mbox{max}\{|\alpha|_{C([0,1]^{2})}, |\beta|_{C([0,1]^{2})}\} \bigl| |\nabla\thetaa_{i}|_{\sigma} \bigr|_{L^{1}(\Omega; \R^N)}
    \nonumber
    \\
    &~ \ds + \frac{2\nu^{2}\max\{|\alpha|_{C([0,1]^{2})},|\beta|_{C([0,1]^{2})}\}}{\delta_{\ast}} \bigr| {\ts\sqrt{\beta(\vv_{i})} \nabla \thetaa_{i}} \bigr|_{L^{2}(\Omega)^{N}}^{2}
    \nonumber
    \\[1ex]
    \le &~ \ds 10^{2}R_{\ast}^{3} \left( \frac{|\alpha|_{C([0,1]^{2})}}{10^{2}R_{\ast}^{2}} + 10^{2}R_{\ast}^{2} + 1 \right) (1 + |\omega_{0}|_{H^{1}(\Omega)}^{2})
    \nonumber
    \\[1ex]
    \le &~ 2 \cdot 10^{4} R_{\ast}^{5} (1 + |\omega_{0}|_{H^{1}(\Omega)}^{2}),
\end{align}
by (A2). Taking the sum of (\ref{inq4-4}) and (\ref{inq4-6}), we can see that 
\begin{equation}\label{inq4-7}
\begin{array}{l}
\ds \frac{1}{2h}(|\vv_{i} - \w_{0}|_{L^{2}(\Omega)^{2}}^{2} - |\vv_{i-1} - \w_{0}|_{L^{2}(\Omega)^{2}}^{2}) + \frac{1}{4}|\nabla \vv_{i}|_{L^{2}(\Omega)^{2N}}^{2} \vspace{3mm}\\
\ds \hspace{5mm} + \int_{\Omega} \gamma(w_{i-1}^{\nu, \sigma}) dx  + \left( \int_{\Omega} \gamma(w_{i}^{\nu, \sigma}) dx - \int_{\Omega} \gamma(w_{i-1}^{\nu, \sigma}) dx \right) \vspace{3mm}\\
\ds \hspace{5mm} + \int_{\Omega}g(\v_{i-1}^{\nu,\sigma})dx + \left( \int_{\Omega}g(\v_{i}^{\nu,\sigma})dx - \int_{\Omega}g(\v_{i-1}^{\nu,\sigma})dx \right)  \vspace{3mm}\\
\ds \hspace{5mm} + \frac{A_{\ast}}{2h} (|\thetaa_{i} - \omega_{0}|_{L^{2}(\Omega)}^{2} - |\thetaa_{i-1} - \omega_{0}|_{L^{2}(\Omega)}^{2}) \vspace{3mm}\\
    \ds \hspace{5mm} + |\alpha|_{C([0,1]^{2})} \bigl(| |\nabla\thetaa_{i}|_{\sigma}|_{L^{1}(\Omega; \R^N)} - | |\nabla \thetaa_{i-1}|_{\sigma} |_{L^{1}(\Omega; \R^N)} \bigr) \vspace{3mm}\\
\ds \hspace{5mm} +  \nu^{2} |\beta|_{C([0,1]^{2})} ( \left| \nabla \thetaa_{i} \right|_{L^{2}(\Omega)^{N}}^{2} - |\nabla \thetaa_{i-1}|_{L^{2}(\Omega)^{N}}^{2})  \vspace{3mm}\\
\ds \hspace{5mm} + \frac{\delta_{\ast
}}{|\alpha|_{C([0,1]^{2})}} \int_{\Omega} \alpha(\vv_{i-1})|\nabla \thetaa_{i-1}|_{\sigma} dx + \nu^{2} \frac{\delta_{\ast
}}{|\beta|_{C([0,1]^{2})}}\int_{\Omega} \beta(\vv_{i-1})|\nabla \thetaa_{i-1}|^{2} dx\vspace{3mm}\\
    \ds \le 4 \cdot 10^{4}R_{\ast}^{5} (1+|\w_{0}|_{H^{1}(\Omega)^{2}}^{2} + |\omega_{0}|_{H^{1}(\Omega)}^{2}) + \frac{c^{2}}{2} |u_{i}|_{L^{2}(\Omega)}^{2}.
\end{array}
\end{equation}
Here, we take
\begin{equation*}
\ds B_{\ast} := \min \left\{ \frac{1}{2},\ \frac{\delta_{\ast
}}{|\alpha|_{C([0,1]^{2})}},\ \frac{\delta_{\ast
    }}{|\beta|_{C([0,1]^{2})}} \right\}\ \ \mbox{ and } \ \ C_{\ast} := 4 \cdot 10^{4} R_{\ast}^{5}, 
\end{equation*}
and the sum of (\ref{inq4-7}) from $i=1$ to $i=m \in \mathbb{N}$, then we get the desired result. 
\hfill $\Box$

\section{Proof of Main Theorem 1}
Let $[\eta_{0},w_{0}, \theta_{0}] \in D_{0}$, the constant $\nu_{\ast}$ obtained in Lemma \ref{lem4-1}, and a fixed constant $\nu_{0} \in [0, \nu_{\ast})$. In this section, we prove Main Theorem 1 through a limiting process for (AP)$_{h}^{\nu, \sigma}$ as $h, \sigma\to 0$ and $\nu \to \nu_{0}$. At first, we recall the auxiliary results for the weighted total variations. 
\begin{lem}[cf. {\cite[Lemma 4.6]{MSW}}]\label{lem5-0}
Let $ \delta_* \in (0, 1) $ be a fixed constant, and let $ I \subset (0, \infty) $ be an open interval. Let $\{\nu_{n}\}_{n=1}^{\infty} \subset (\nu_{0},\nu_{\ast})$, $\{\sigma_{n}\}_{n=1}^{\infty} \subset (0,1)$ with $\nu_{n} \downarrow \nu_{0}$, $\sigma_{n} \downarrow 0$ as $n \to \infty$, respectively. Also, let us assume that 
\begin{equation*}
\left \{ \begin{array}{l}
\ds \varrho \in C(\overline{I};L^{2}(\Omega)) \cap L^{\infty}(I;H^{1}(\Omega)) \cap L^{\infty}(I \times \Omega), \{\varrho_{n}\}_{n=1}^{\infty} \subset L^{2}(I;L^{2}(\Omega)), \vspace{3mm}\\
\ds \varrho \ge 0 \mbox{ and } \rho_{n} \ge 0, \mbox{ a.e. in } I \times \Omega, \mbox{ for all } n \in \N, \vspace{3mm}\\
\ds \zeta \in C(\overline{I}; L^{2}(\Omega)), \{\zeta_{n}\}_{n=1}^{\infty} \subset L^{2}(I;H^{1}(\Omega)), \vspace{3mm}\\
\ds \varrho_{n}(t) \to \varrho(t) \mbox{ in } L^{2}(\Omega) \mbox{ and weakly in } H^{1}(\Omega), \mbox{ for a.e. } t \in I, \mbox{ as } n \to \infty, \vspace{3mm}\\
\ds \zeta_{n}(t) \to \zeta(t) \mbox{ in } L^{2}(\Omega), \mbox{ for a.e. } t \in I, \mbox{ as } n \to \infty. 
\end{array}
\right. 
\end{equation*}
In addition, let us assume that 
\begin{equation*}
\varrho \ge \delta_{\ast} \mbox{ a.e. in } I \times \Omega, \mbox{ or } L_{0}:= \sup_{n \in \N}|\nabla \zeta_{n}|_{L^{1}(I;L^{2}(\omega ; \R^{N}))} < \infty. 
\end{equation*}
Then, 
\begin{equation*}
\liminf_{n \to \infty} \Phi_{\nu_{n}, \sigma_{n}}^{I}(\varrho_{n}; \zeta_{n}) \ge \liminf_{n \to \infty} \Phi_{\nu_{0}}^{I}(\varrho_{n}; \zeta_{n}) \ge \Phi_{\nu_{0}}^{I}(\varrho; \zeta). 
\end{equation*}
\end{lem}

\begin{lem}[cf. {\cite[Lemma 4.9]{MSW}} and {\cite[Lemma 6.1]{SWY14}}]\label{lem5-1}
Let $ \delta_* \in (0, 1) $ be a fixed constant, and let $ I \subset (0, \infty) $ be an open interval. Assume that
\begin{equation*}
\left\{ \hspace{-5ex} \parbox{15cm}{
\vspace{-2ex}
\begin{itemize}
\item[]$ \varrho \in C(\overline{I}; L^2(\Omega)) \cap L^\infty(I; H^1(\Omega)) \cap L^\infty(I \times \Omega) $, $ \{ \varrho_n \, | \, n \in \N \} \subset L^2(I; L^2(\Omega)) $,
\vspace{-1ex}
\item[]$ \varrho \geq \delta_* $ and $ \varrho_n \geq \delta_* $, a.e.\ in $ I \times \Omega $, for all $ n \in \N $,
\vspace{-1ex}
\item[]$ \varrho_n(t) \to \varrho(t) $ in $ L^2(\Omega) $ and weakly in $ H^1(\Omega) $, as $ n \to \infty $, a.e.\ $ t \in I $,
\vspace{-2ex}
\end{itemize}
} \right.
\end{equation*}
and
\begin{equation*}
\left\{ \hspace{-5ex} \parbox{14cm}{
\vspace{-2ex}
\begin{itemize}
\item[]$ \zeta \in C(\overline{I}; L^2(\Omega)) \cap L^1(I; BV(\Omega)) $, $ \{ \zeta_n \, | \, n \in \N \} \subset L^2(I; H^1(\Omega)) $,
\vspace{-1ex}
\item[]$ \zeta_n(t) \to \zeta(t) $ in $ L^2(\Omega) $ as $ n \to \infty $, \ a.e.\ $ t \in I $.
\vspace{-2ex}
\end{itemize}
} \right.
\end{equation*}
Then the functions
\begin{equation*}
t \in I \mapsto \int_\Omega d[\varrho(t)|D\zeta(t)|], \mbox{ and } 
t \in I \mapsto \int_\Omega \varrho_n(t)|\nabla \zeta_n(t)| \, dx, ~ n \in \N,
\end{equation*}
are integrable. Moreover, if
\begin{equation*}
\int_I \int_\Omega \varrho_n(t) |\nabla \zeta_n(t)|_{\sigma_{n}} \, dx \, dt \to \int_I \int_\Omega d[\varrho(t) |D \zeta(t)|]dt
\end{equation*}
as $ n \to \infty $, and
\begin{equation*}
\left\{ \hspace{-5ex} \parbox{15cm}{
\vspace{-2ex}
\begin{itemize}
\item[]$ \omega \in C(\overline{I}; L^2(\Omega)) \cap L^\infty(I; H^1(\Omega)) \cap L^\infty(I \times \Omega) $ and $ \{ \omega_n \, | \, n \in \N \} \subset L^2(I; L^2(\Omega)) $,
\vspace{-4ex}
\item[]$ \{ \omega_n \, | \, n \in \N \} $ is a bounded sequence in $ L^\infty(I \times \Omega) $,
\vspace{-1ex}
\item[]$ \omega_n(t) \to \omega(t) $ in $ L^2(\Omega) $ and weakly in $ H^1(\Omega) $ as $ n \to \infty $, a.e.\ $ t \in I $,
\vspace{-2ex}
\end{itemize}
} \right.
\end{equation*}
then
\begin{equation*}
\int_I \int_\Omega \omega_n(t) |\nabla \zeta_n(t)|_{\sigma_{n}} \, dx \, dt \to \int_I \int_\Omega d [\omega(t) |D \zeta(t)|]dt
\end{equation*}
as $ n \to \infty $.
\end{lem}
\bigskip

Let $\nu_{\ast}$, $h_{\ast}:=h_{1}^{\dagger} \in (0,1)$ be the constants as in (\ref{nu-ast}) and (\ref{h1}). We take $\{[\tilde{\v}_{0}^{\sigma}, \tilde{\theta}_{0}^{\sigma}]\}_{\sigma \in (0,1)} = \{[\tilde{w}_{0}^{\sigma}, \tilde{\eta}_{0}^{\sigma}, \tilde{\theta}_{0}^{\sigma}]\}_{\sigma \in (0,1)}$ such that 
\begin{equation}\label{5-1}
\left\{ \begin{array}{l}
\ds [\tilde{\v}_{0}^{\sigma}, \tilde{\theta}_{0}^{\sigma}] \in D_{1}^{\ast}(\theta_{0}) \ \ \mbox{ for all } \sigma \in (0,1), \vspace{3mm}\\
\ds [\tilde{\v}_{0}^{\sigma}, \tilde{\theta}_{0}^{\sigma}] \to [\v_{0}, \theta_{0}]\ \ \mbox{ in } L^{2}(\Omega)^{3}\ \mbox{ as } \sigma \to 0.
\end{array}\right. 
\end{equation}
Then, for any $\nu \in (\nu_{0},\nu_{\ast})$, $\sigma \in (0,1)$, and $h \in (0,h_{\ast})$,  Lemma \ref{thm1} guarantees the existence of a unique solution $\{[ \tilde{\v}_{i}^{\nu, \sigma}, \tilde{\theta}_{i}^{\nu, \sigma}]\}_{i=1}^{\infty}$ to (AP)$_{h}^{\nu, \sigma}$ in the case of $[\v_{0}^{\nu, \sigma}, \theta_{0}^{\nu,\sigma}] = [\tilde{\v}_{0}^{\sigma}, \tilde{\theta}_{0}^{\sigma}]$. Also, we take three kinds of time interpolations $ [\overline{\v}_h^{\nu, \sigma}, \overline{\theta}_h^{\nu, \sigma}] \in L_{\rm loc}^\infty([0, \infty); H^1(\Omega))^{3} $, $ [\underline{\v}_h^{\nu, \sigma}, \underline{\theta}_h^{\nu, \sigma}] \in L_{\rm loc}^\infty([0, \infty); H^1(\Omega))^{3} $, and $ [\widehat{\v}_h^{\nu, \sigma}, \widehat{\theta}_h^{\nu, \sigma}] \in W_{\rm loc}^{1, \infty}([0, \infty); H^1(\Omega))^{3} $, by letting
\begin{equation}\label{5-2}
\left \{
\begin{array}{ll}
\ds [\overline{\v}_{h}^{\nu, \sigma}(t), \overline{\theta}_{h}^{\nu, \sigma}(t)] := [\tilde{\v}_{i}^{\nu, \sigma}, \tilde{\theta}_{i}^{\nu, \sigma}], & \mbox{if $ t \in ((i -1)h, ih] \cap [0, \infty) $ with some $ i \in \mathbb{Z} $,}
\\[2ex]
[\underline{\v}_{h}^{\nu, \sigma}(t), \underline{\theta}_{h}^{\nu, \sigma}(t)] := [\tilde{\v}_{i-1}^{\nu, \sigma}, \tilde{\theta}_{i-1}^{\nu, \sigma}], & \mbox{if $ t \in [(i -1)h, ih) $ with some $ i \in \N $,}
\\[1ex]
\multicolumn{2}{l}{
\ds [\widehat{\v}_{h}^{\nu, \sigma}(t), \widehat{\theta}_{h}^{\nu, \sigma}(t)] := \frac{ih -t}{h}[\tilde{\v}_{i-1}^{\nu, \sigma}, \tilde{\theta}_{i-1}^{\nu, \sigma}] + \frac{t-(i -1)h}{h}[\tilde{\v}_{i}^{\nu, \sigma}, \tilde{\theta}_{i}^{\nu, \sigma}],
}
\\[1ex]
& \mbox{if $ t \in [(i -1)h, ih) $ with some $ i \in \N $,}
\end{array}
\right.
\end{equation}
for all $ t \geq 0 $. Moreover, we define $\overline{u}_{h} \in L^{2}_\mathrm{loc}([0,\infty);L^{2}(\Omega))$ by 
\begin{equation*}
\overline{u}_{h} := u_{i}
\mbox{ if }
t \in ((i-1)h,ih] \cap [0,\infty)
\mbox{ with some $0 \le i \in \mathbb{Z}$}. 
\end{equation*}
Then, due to Lemma \ref{thm1}, (\ref{5-1}) implies that
\begin{equation}\label{5-4}
\begin{array}{c}
\left\{
\begin{array}{l}
\ds \overline{\v}_h^{\nu, \sigma}(t),\ \underline{\v}_h^{\nu, \sigma}(t),\ \widehat{\v}_h^{\nu, \sigma}(t) \in [0,1]^{2},
\\[0.5ex]
\ds \max \left\{ |\overline{\theta}_h^{\nu, \sigma}(t)|, ~ |\underline{\theta}_h^{\nu, \sigma}(t)|, ~ |\widehat{\theta}_h^{\nu, \sigma}(t)| \right\} \leq |\tilde{\theta}_{0}^{\sigma}|_{L^\infty(\Omega)} \leq |\theta_0|_{L^\infty(\Omega)},
\end{array}
\right. 
\ \\
\\[-2.5ex]
\end{array}
\end{equation}
a.e. in $ \Omega $, for all $ t \geq 0 $, $ \nu \in (\nu_{0}, \nu_*) $, $\sigma \in (0,1)$, and $ h \in (0, h_*) $. By using these interpolations, the inequality (\ref{AP_diss}) of energy dissipation leads to 
\begin{equation}\label{5-5}
\begin{array}{l}
    \ds \frac{1}{2} \int_{s}^{t} |(\widehat{\v}_{h}^{\nu, \sigma})_{t}(\tau)|_{L^{2}(\Omega)^{2}}^{2}d \tau  + \int_{s}^{t} \bigl|{\ts \sqrt{\alpha_{0}(\overline{\v}_{h}^{\nu, \sigma}(\tau))}(\widehat{\theta}_{h}^{\nu, \sigma})_{t}(\tau) } \bigr|_{L^{2}(\Omega)}^{2} d\tau \vspace{3mm}\\
\ds   + \F_{\nu, \sigma}(\overline{\v}_{h}^{\nu, \sigma}(t), \overline{\theta}_{h}^{\nu, \sigma}(t)) + c \int_{\Omega} {u}^\dag \overline{w}_{h}^{\nu, \sigma}(t) dx - c^{2}\int_{0}^{t} \int_{\Omega}|\overline{u}_{h}(\tau)-{u}^\dag|^{2}dx d\tau\vspace{3mm}\\
\ds \le \F_{\nu, \sigma}(\underline{\v}_{h}^{\nu, \sigma}(s), \underline{\theta}_{h}^{\nu, \sigma}(s)) + c \int_{\Omega} {u}^\dag \underline{w}_{h}^{\nu, \sigma}(s) dx - c^{2}\int_{0}^{s} \int_{\Omega}|\overline{u}_{h}(\tau)-{u}^\dag|^{2}dx d\tau, 
\end{array}
\end{equation}
for any $ 0 \le s \le t < +\infty $. Similarly, (\ref{AP_diss2}) and (\ref{lem1-inq}) derive 
\begin{equation}\label{5-6}
\begin{array}{c}
    \ds \frac{1}{2} \int_{0}^{t} \tau |(\widehat{\v}_{h}^{\nu, \sigma})_{t}(\tau)|_{L^{2}(\Omega)^{2}}^{2}d \tau  + \int_{0}^{t} \tau \bigl| {\ts \sqrt{\alpha_{0}(\overline{\v}_{h}^{\nu, \sigma}(\tau))}(\widehat{\theta}_{h}^{\nu, \sigma})_{t}(\tau)} \bigr|_{L^{2}(\Omega)}^{2} d\tau \vspace{3mm}\\
\ds \hspace{0mm} + t \F_{\nu, \sigma}(\overline{\v}_{h}^{\nu, \sigma}(t), \overline{\theta}_{h}^{\nu, \sigma}(t)) + ct({u}^\dag, \overline{w}_{h}^{\nu,\sigma}(t))_{L^{2}(\Omega)} 
\ds 
    \ds \le \int_{0}^{t+h} \F_{\nu, \sigma}(\underline{\v}_{h}^{\nu, \sigma}(\tau), \underline{\theta}_{h}^{\nu, \sigma}(\tau))d\tau \vspace{3mm}\\
    \ds + c \int_{0}^{t} ({u}^\dag, \underline{w}_{h}^{\nu,\sigma}(\tau))_{L^{2}(\Omega)} d\tau +c^{2} (t+h) \int_{0}^{t+h} 
|\overline{u}_{h}(\tau) - {u}^\dag|_{L^{2}(\Omega)}^{2}
    d\tau + 2ch|u^{\dag}|\mathscr{L}^{N}(\Omega), 
\end{array}
\end{equation}
and 
\begin{equation}\label{5-7}
\begin{array}{l}
\ds \frac{1}{2}(|\overline{\v}_{h}^{\nu, \sigma}(t) -\w_{0}|_{L^{2}(\Omega)^{2}}^{2} + A_{\ast} |\overline{\theta}_{h}^{\nu, \sigma}(t) -\omega_{0}|_{L^{2}(\Omega)}^{2}) + {\frac{B_*}{2}} \int_{0}^{t} \mathscr{F}_{\nu, \sigma}(\underline{\v}_{h}^{\nu, \sigma}(\tau), \underline{\theta}_{h}^{\nu, \sigma}(\tau))d\tau \vspace{3mm}\\
\ds \le \ds \frac{1}{2}(|\tilde{\v}_{0}^{\sigma}-\w_{0}|_{L^{2}(\Omega)^{2}}^{2} + A_{\ast} |\tilde{\theta}_{0}^{\sigma}-\omega_{0}|_{L^{2}(\Omega)}^{2}) +\frac{h}{B_*} \mathscr{F}_{\nu, \sigma}(\tilde{\v}_{0}^{\sigma}, \tilde{\theta}_{0}^{\sigma})\vspace{3mm}\\
    \ds\hspace{5mm} + 2t C_*(1 +|\w_0|_{H^1(\Omega)^{2}}^2 +|\omega_0|_{H^1(\Omega)}^2) + \frac{c^{2}}{2} \int_{0}^{t+h} |\overline{u}_{h}(\tau)|_{L^{2}(\Omega)}^{2}d\tau,
\end{array}
\end{equation}
for all $0 \le t < +\infty$, $\nu \in (\nu_{0},\nu_{\ast})$, $\sigma \in (0,1)$, and $h \in (0,h_{\ast})$, respectively. 

Using (\ref{5-1}) and a diagonal argument, we can obtain sequences $\{ \nu_{n} \}_{n=1}^{\infty} \subset (\nu_{0},\nu_{\ast})$, $\{\sigma_{n}\}_{n=1}^{\infty} \subset (0,1)$, and $\{h_{n}\}_{n=1}^{\infty} \subset (0, h_{\ast})$, such that 
\begin{enumerate}
\item[(1)] if $\nu_{0} > 0$, then 
\begin{equation}\label{5-7.4}
\left\{ \begin{array}{l}
\ds \nu_{n} := \nu_{0},\ \ 0 < \sigma_{n+1} < \sigma_{n} < 2^{-n},\ \ 0 < h_{n+1} < h_{n} < h_{\ast} 2^{-n}, \vspace{3mm}\\
\ds 0 \le h_{n} \F_{\nu_{n}, \sigma_{n}}(\tilde{\v}_{0}^{\sigma_{n}} , \tilde{\theta}_{0}^{\sigma_{n}}) < 2^{-n}, 
\end{array} \right. 
\end{equation}
for all $n \in \N$; 
\item[(2)] if $\nu_{0} = 0$, then 
\begin{equation}\label{5-7.5}
\left\{ \begin{array}{l}
\ds 0 = \nu_{0} < \nu_{n+1} < \nu_{n} < \nu_{\ast}2^{-n},\ \ \sigma_{n} := \nu_{n},\ \ 0 < h_{n+1} < h_{n} < h_{\ast} 2^{-n}, \vspace{3mm}\\
\ds 0 \le h_{n} \F_{\nu_{n}, \sigma_{n}}(\tilde{\v}_{0}^{\sigma_{n}} , \tilde{\theta}_{0}^{\sigma_{n}}) < 2^{-n}, 
\end{array} \right. 
\end{equation}
for all $n \in \N$. 
\end{enumerate}
According to (\ref{5-4})-(\ref{5-7.5}), the sequences 
\begin{equation}\label{5-7.51}
\left\{ 
\begin{array}{l}
    \ds \{ [\overline{\v}_{n}(t), \overline{\theta}_{n}(t)] \}_{n=1}^{\infty} := \{ [\overline{\v}_{h_{n}}^{\nu_{n}, \sigma_{n}}(t), \overline{\theta}_{h_{n}}^{\nu_{n}, \sigma_{n}}(t)] \}_{n=1}^{\infty} \subset L_\mathrm{loc}^\infty([0, \infty); H^1(\Omega)^3),
    \\[1.5ex]
\ds \{ [\underline{\v}_{n}(t), \underline{\theta}_{n}(t)] \}_{n=1}^{\infty} := \{ [\underline{\v}_{h_{n}}^{\nu_{n}, \sigma_{n}}(t), \underline{\theta}_{h_{n}}^{\nu_{n}, \sigma_{n}}(t)] \}_{n=1}^{\infty}\subset L_\mathrm{loc}^\infty([0, \infty); H^1(\Omega)^3),
    \\[1.5ex]
    \ds \{ [\widehat{\v}_{n}(t), \widehat{\theta}_{n}(t)] \}_{n=1}^{\infty} := \{ [\widehat{\v}_{h_{n}}^{\nu_{n}, \sigma_{n}}(t), \widehat{\theta}_{h_{n}}^{\nu_{n}, \sigma_{n}}(t)] \}_{n=1}^{\infty} \subset W_\mathrm{loc}^{1, \infty}([0, \infty); H^1(\Omega)^3),
    \\[1.5ex]
    \ds \{ [\v_{0,n}, \theta_{0,n}] \}_{n=1}^{\infty} := \{ [\tilde{\v}_{0}^{\nu_n, \sigma_{n}}, \tilde{\theta}_{0}^{\nu_n, \sigma_{n}}] \}_{n=1}^{\infty} \subset D_1(\theta_0),
\end{array}
\right. 
\end{equation}
\begin{equation}\label{5-7.51-01}
    \{ \overline{u}_n \}_{n = 1}^\infty := \{ \overline{u}_{h_n} \}_{n = 1}^\infty \subset L_\mathrm{loc}^2([0, \infty); L^2(\Omega)),
\end{equation}
and sequences $ \{ \overline{\F}_n^{{u}^\dag} \}_{n = 1}^\infty, \{ \underline{\F}_n^{{u}^\dag} \}_{n = 1}^\infty \subset L_\mathrm{loc}^1([0, \infty)) $, defined as:
\begin{equation}\label{5-7.52}
    \begin{array}{c}
\left\{ 
\begin{array}{l}
    \ds \overline{\F}_n^{{u}^\dag}(t):= \F_{\nu_{n},\sigma_{n}}(\overline{\v}_{n}(t), \overline{\theta}_{n}(t)) 
    \\[1ex]
    \ds \qquad + c \int_{\Omega} {u}^\dag \overline{w}_{n}(t) dx
 - c^{2}\int_{0}^{t} \int_{\Omega}|\overline{u}_{h_{n}}(\tau)-{u}^\dag|^{2}dx d\tau, 
    \\[2ex]
    \ds \underline{\F}_n^{{u}^\dag}(t) := \F_{\nu_{n},\sigma_{n}}(\underline{\v}_{n}(t), \underline{\theta}_{n}(t)) 
    \\[1ex]
    \ds \qquad + c \int_{\Omega} {u}^\dag \underline{w}_{n}(t) dx - c^{2}\int_{0}^{t} \int_{\Omega}|\overline{u}_{h_{n}}(\tau)-{u}^\dag|^{2}dx d\tau,
\end{array}
\right. 
        \ \\[-1ex]
        \ \\
        \mbox{for all $ t \geq 0 $, and $ n = 1, 2, 3, \dots $,}
    \end{array}
\end{equation}
satisfy the following properties: 
\begin{description}
    \item[($\sharp$a)] $\{ [\overline{\v}_{n}(t), \overline{\theta}_{n}(t)],\ [\underline{\v}_{n}(t), \underline{\theta}_{n}(t)],\  [\widehat{\v}_{n}(t), \widehat{\theta}_{n}(t)] \}_{n=1}^{\infty} \subset D_{1}(\theta_{0})$, for all $t \ge 0$; 
\item[($\sharp$b)] $\{ [\overline{\v}_{n}, \overline{\theta}_{n}]\}_{n=1}^{\infty}$ and $\{[\underline{\v}_{n}, \underline{\theta}_{n}] \}_{n=1}^{\infty} $ are bounded in $L^{\infty}_\mathrm{loc}((0,\infty);L^{2}(\Omega)^{3}) \cap L^{\infty}(Q)^{3}$, 
$\{ [\widehat{\v}_{n}, \widehat{\theta}_{n}] \}_{n=1}^{\infty}$ is bounded in $C_\mathrm{loc}((0,\infty);L^{2}(\Omega)^{3}) \cap W^{1,2}_\mathrm{loc}((0,\infty);L^{2}(\Omega)^{3}) \cap L^{\infty}_\mathrm{loc}((0,\infty) ; H^{1}(\Omega) \times H^{1}(\Omega) \times BV(\Omega)) \cap L^{\infty}(Q)^{3}$. Also, 
$\{\nu_{n}\overline{\theta}_{n}\}_{n=1}^{\infty}$, $\{\nu_{n}\underline{\theta}_{n}\}_{n=1}^{\infty}$, and $\{\nu_{n}\widehat{\theta}_{n}\}_{n=1}^{\infty}$ are bounded in $L^{\infty}_\mathrm{loc}((0,\infty) ; H^{1}(\Omega))$;
\item[($\sharp$c)] $\{ \underline{\F}_n^{{u}^\dag}\}_{n=1}^{\infty}$ and $\{ \underline{\F}_n^{{u}^\dag}\}_{n=1}^{\infty}$ are sequences of nonincreasing functions on $(0, \infty)$, which are bounded in $L^{1}_\mathrm{loc}([0,\infty))$ and $BV_\mathrm{loc}((0,\infty))$;
\item[($\sharp$d)] $h_{n} \F_{\nu_{n}, \sigma_{n}}(\v_{0,n}, \theta_{0,n}) \to 0$ as $n \to \infty$;
    
    \item[($\sharp$e)] $ \overline{u}_n \to u $ in $ L_\mathrm{loc}^2([0, \infty); L^2(\Omega)) $ as $ n \to \infty $.
\end{description}
Taking into account the compactness results as in \cite[Chapter 3]{AFP} and \cite[Corollary 4]{Simon}, there exists $[\v, \theta] \in C_\mathrm{loc}((0,\infty);L^{2}(\Omega)^{3})$ and subsequences (not relabeled) of $\{[\widehat{\v}_{n},\widehat{\theta}_{n}]\}_{n=1}^{\infty}$ such that 
\begin{equation}\label{5-8}
\left\{
\begin{array}{l}
\ds \v \in W^{1,2}_\mathrm{loc}((0,\infty);L^{2}(\Omega)^{2}) \cap L^{\infty}_\mathrm{loc}((0,\infty);H^{1}(\Omega)^{2}), \vspace{3mm}\\
\ds \theta \in W^{1,2}_\mathrm{loc}((0,\infty) ; L^{2}(\Omega)),\ |D\theta(\cdot)|(\Omega) \in L^{\infty}_\mathrm{loc}((0,\infty)), \vspace{3mm}\\
\ds \nu_{0}\theta \in L_\mathrm{loc}^{\infty}((0,\infty);H^{1}(\Omega)), \vspace{3mm}\\
\ds 0 \le w \le 1,\ 0\le \eta \le 1, \mbox{ and } |\theta| \le |\theta_{0}|_{L^{\infty}(\Omega)}\ \mbox{ a.e. in } Q, 
\end{array}
\right. 
\end{equation}
and 
\begin{equation}\label{5-9}
\left\{ 
\begin{array}{l}
\ds \widehat{\v}_{n} \to \v \mbox{ in } C_\mathrm{loc}((0,\infty);L^{2}(\Omega)^{2}), \mbox{ weakly in } W^{1,2}_\mathrm{loc}((0,\infty);L^{2}(\Omega)^{2}), \vspace{3mm}\\
\ds \hspace{13mm} \mbox{ weakly-}\ast \mbox{ in } L^{\infty}_\mathrm{loc}((0,\infty);H^{1}(\Omega)^{2}), \ds \mbox{ and weakly-}\ast \mbox{ in } L^{\infty}(Q)^{2}, \vspace{3mm}\\
\ds \widehat{\v}_{n}(t) \to \v(t) \mbox{ in } L^{2}(\Omega)^{2}, \mbox{ weakly in } H^{1}(\Omega)^{2} \mbox{ for any } t > 0, \vspace{3mm}\\
\ds \widehat{\theta}_{n} \to \theta \mbox{ in } C_\mathrm{loc}((0,\infty) ; L^{2}(\Omega)), \mbox{ weakly in } W^{1,2}_\mathrm{loc}((0,\infty);L^{2}(\Omega)), \vspace{3mm}\\
\ds \hspace{13mm} \mbox{ and weakly-}\ast \mbox{ in } L^{\infty}(Q), \vspace{3mm}\\
\ds \widehat{\theta}_{n}(t) \to \theta(t) \mbox{ in } L^{2}(\Omega), \mbox{ and weakly-}\ast \mbox{ in } BV(\Omega), \mbox{ for any } t>0,
\end{array}
\right. 
\end{equation}
as $n \to \infty$. Moreover, by the energy inequality (\ref{5-5}), we obtain
\begin{equation}\label{5-10}
\left\{ 
\begin{array}{l}
\ds \max \{ |\overline{\v}_{n}(t) - \widehat{\v}_{n}(t)|_{L^{2}(\Omega)^{2}},\ |\underline{\v}_{n}(t) - \widehat{\v}_{n}(t)|_{L^{2}(\Omega)^{2}} \} \le \int_{(i-1)h}^{ih} |(\widehat{\v}_{n})_{t}(t)|_{L^{2}(\Omega)^{2}} dt, 
    \\[2ex]
\ds \max \{ |\overline{\theta}_{n}(t) - \widehat{\theta}_{n}(t)|_{L^{2}(\Omega)},\ |\underline{\theta}_{n}(t) - \widehat{\theta}_{n}(t)|_{L^{2}(\Omega)} \} \le \int_{(i-1)h}^{ih} |(\widehat{\theta}_{n})_{t}(t)|_{L^{2}(\Omega)} dt,
\end{array}
\right.
\end{equation}
for $ n = 1, 2, 3, \dots $. 
Therefore, the following convergences
\begin{equation}\label{5-11}
\left\{ 
\begin{array}{l}
\ds \overline{\v}_{n} \to \v,\ \underline{\v}_{n} \to \v \mbox{ in } L^{\infty}_\mathrm{loc}((0,\infty);L^{2}(\Omega)^{2}),
    \\[1ex]
\ds \qquad \mbox{weakly-}\ast \mbox{ in } L^{\infty}_\mathrm{loc}((0,\infty);H^{1}(\Omega)^{2}), 
    \\
    \qquad \mbox{and weakly-}\ast \mbox{ in } L^{\infty}(Q)^{2}, 
    \\[1ex]
\ds \overline{\v}_{n}(t) \to \v(t),\  \underline{\v}_{n}(t) \to \v(t) \mbox{ in } L^{2}(\Omega)^{2}, 
    \\
    \qquad \mbox{weakly in } H^{1}(\Omega)^{2}, \mbox{ for any } t>0, 
\end{array}
\right.
\end{equation}
and 
\begin{equation}\label{5-12}
\left\{ 
\begin{array}{l}
\ds \overline{\theta}_{n} \to \theta,\  \underline{\theta}_{n} \to \theta \mbox{ in } L^{\infty}_\mathrm{loc}((0,\infty);L^{2}(\Omega)),
    \\
    \qquad \mbox{and weakly-}\ast \mbox{ in } L^{\infty}(Q), 
    \\[1ex]
\ds \overline{\theta}_{n}(t) \to \theta(t),\ \underline{\theta}_{n}(t) \to \theta(t) \mbox{ in } L^{2}(\Omega), 
    \\
    \qquad \mbox{and weakly-}\ast \mbox{ in } BV(\Omega), \mbox{ for any } t>0,
    \\[1ex]
    \bigl({\ts \sqrt{\beta(\underline{\v}_{n})}}\nabla (\nu_n \underline{\theta}_{n})\bigr)(t) \to \bigl({\ts \sqrt{\beta(\v)} \nabla(\nu_0 \theta)} \bigr)(t)
    \\
    \qquad \mbox{weakly in $ L^2(\Omega)^N $, for any $ t > 0 $,}
\end{array}
\right. 
\end{equation}
hold as $n \to \infty$. If $\nu_{0} > 0$, then additional convergences follow that
\begin{equation}\label{5-12.1}
\left\{ \begin{array}{l}
\ds \widehat{\theta}_{n},\ \overline{\theta}_{n},\ \underline{\theta}_{n} \to \theta \mbox{ weakly-}\ast \mbox{ in } L^{\infty}_\mathrm{loc}((0,\infty); H^{1}(\Omega)), 
    \\[1ex] 
\ds \widehat{\theta}_{n}(t),\ \overline{\theta}_{n}(t),\ \underline{\theta}_{n}(t) \to \theta(t) \mbox{ weakly in } H^{1}(\Omega) \mbox{ for any } t> 0, 
    \\[1ex]
    \ds (\alpha(\underline{\v}_{n})\nabla\underline{\theta}_{n})(t) \to (\alpha(\v)\nabla\theta)(t) \mbox{ and } ({\ts \sqrt{\beta(\underline{\v}_{n})}}\nabla \underline{\theta}_{n})(t) \to ({\ts \sqrt{\beta(\v)}\nabla\theta})(t)
    \\[1ex]
\ds \hspace{5mm} \mbox{ weakly in } L^{2}(\Omega)^{N}, \mbox{ for any } t > 0,
\end{array} \right. 
\end{equation}
as $n \to \infty$. Moreover, there exists a function $\mathscr{J}^{{u}^\dag}_{\ast} \in BV_\mathrm{loc}((0,\infty))$ such that
\begin{align}\label{5-12.3}
    \ds \underline{\F}_{n}^{{u}^\dag} \to \mathscr{J}^{{u}^\dag}_{\ast} &~ \mbox{weakly-}\ast \mbox{ in } BV_\mathrm{loc}((0,\infty)),  
    \nonumber
    \\
    &~ \ds \mbox{weakly-}\ast \mbox{ in } L^{\infty}_\mathrm{loc}((0,\infty)), \mbox{ and a.e. in } (0,\infty), 
\end{align}
as $n \to \infty$, by taking a suitable subsequence if necessary.
Thus, we get the convergence results for the approximate sequences.
\medskip

Next, we prove that the limit function $[\v,\theta]$ satisfies the variational inequalities in (S2) and (S3). To see this, we define the time-dependent weighted total variation (cf. \cite{MS14}), and refer the corresponding convergence results (cf. \cite[Theorem 4.8]{MSW}, \cite[Lemma 5.1, Lemma 6.3]{SWY14}, and \cite[Main Theorem 1]{SW15}). 
\begin{lem}\label{lem-gamma}
Let $I \subset (0,\infty)$ be a fixed bounded open interval, and let $\Phi_{\nu_{0}}^{I} : L^{2}(I;L^{2}(\Omega)) \to [0,\infty]$ and $\Phi_{\nu_{n}, \sigma_{n}}^{I} : L^{2}(I;L^{2}(\Omega)) \to [0,\infty]$ be functionals defined as
\begin{equation*}\label{5-12.4}
\zeta \in L^{2}(I;L^{2}(\Omega)) \mapsto \ds \Phi_{\nu_{0}}^{I}(\v;\zeta) := \int_{I}\Phi_{\nu_{0}}(\v(t);\zeta(t))dt \in [0, \infty], 
\end{equation*}
and 
\begin{equation*}
    \zeta \in L^{2}(I;L^{2}(\Omega)) \mapsto \Phi_{\nu_{n}, \sigma_{n}}^{I}(\overline{\v}_{n};\zeta) := \int_{I} \Phi_{\nu_{n}}^{\sigma_{n}}(\overline{\v}_{n}(t);\zeta(t))dt \in [0, \infty],
\end{equation*}
    for $\v = [w,\eta] \in L^{\infty}(I;H^{1}(\Omega)^{2}) \cap L^{\infty}(I \times \Omega)^{2}$ and $\overline{\v}_{n} = [\overline{w}_{n},\overline{\eta}_{n}] \in L^{\infty}(I;H^{1}(\Omega)^{2}) \cap L^{\infty}(I \times \Omega)^{2}$, $ n \in \N $, as in (\ref{5-8})-(\ref{5-12}). Then, the following two statements hold: 
\begin{description}
\item[(G-1)] $\Phi_{\nu_{0}}^{I}(\v;{}\cdot\,)$, and $\Phi_{\nu_{n}, \sigma_{n}}^{I}(\overline{\v}_{n};{}\cdot\,)$, $n \in \N$, are proper l.s.c. and convex functions on $L^{2}(I; L^{2}(\Omega))$ such that 
\begin{itemize}
\item if $\nu_{0} = 0$, then $\D(\Phi_{\nu_{0}}^{I}(\v;{}\cdot\,)) = L^{1}(I ; BV(\Omega)) \cap L^{2}(I ; L^{2}(\Omega))$, 
\item if $\nu_{0} > 0$, then $\D(\Phi_{\nu_{0}}^{I}(\v;{}\cdot\,)) = \D(\Phi_{\nu_{n},\sigma_{n}}^{I}(\overline{\v}_{n};{}\cdot\,)) = L^{2}(I ; H^{1}(\Omega))$, for all $n \in \N$. 
\end{itemize}
\item[(G-2)] the sequence $\{\Phi_{\nu_{n}, \sigma_{n}}^{I}(\overline{\v}_{n};{}\cdot\,)\}_{n=1}^{\infty}$ converges to $\Phi_{\nu_{0}}^{I}(\v;{}\cdot\,)$ on $L^{2}(I;L^{2}(\Omega))$, in the sense of $\Gamma$-convergence, as $n \to \infty$. 
\end{description}
\end{lem}
\begin{rem}
    \begin{em}
If $\nu_{0}>0$, then the sequence $\{\Phi_{\nu_{n}, \sigma_{n}}^{I}(\overline{\v}_{n};{}\cdot\,)\}_{n=1}^{\infty}$ converges to $\Phi_{\nu_{0}}^{I}(\v;{}\cdot\,)$ on $L^{2}(I;L^{2}(\Omega))$, in the sense of Mosco-convergence, as $n \to \infty$ (cf. \cite[Lemma 5.1]{SWY14}). 
    \end{em}
\end{rem}

Let $I$ be any bounded open interval such that $I \subset \subset (0,\infty)$. By (\ref{apEq.2}) and (\ref{apEq.3}), the sequences, as in (\ref{5-9})-(\ref{5-12}), satisfy  the following two variational inequalities: 
\begin{equation}\label{5-12.5}
\begin{array}{l}
    \ds \int_{I} ((\widehat{\v}_{n})_{t}(t), \overline{\v}_{n}(t) - \bm{\varpi}(t))_{L^{2}(\Omega)^{2}} dt + \int_{I} (\nabla \overline{\v}_{n}(t), \nabla (\overline{\v}_{n} - \bm{\varpi})(t))_{L^{2}(\Omega)^{2}} dt \v\\[2ex]
\ds \quad + \int_{I} ([\nabla G](\overline{u}_{h_{n}};\overline{\v}_{n})(t), \overline{\v}_{n}(t) - \bm{\varpi}(t))_{L^{2}(\Omega)^{2}} dt + \int_{I} \int_{\Omega} \gamma (\overline{\v}_{n}(t)) dxdt 
    \\[2ex]
\ds \quad + \int_{I} \int_{\Omega} \bigl( \overline{\v}_{n} - \bm{\varpi})(t) \cdot (|\nabla\underline{\theta}_{n}(t)|_{\sigma_{n}}[\nabla\alpha](\overline{\v}_{n}(t)) + \nu_{n}^{2}|\nabla\underline{\theta}_{n}(t)|^{2} [\nabla\beta](\overline{\v}_{n}(t)) \bigr) dxdt 
    \\[2ex]
\ds \le \int_{I}\int_{\Omega} \gamma(\bm{\varpi}(t)) dxdt
\end{array}
\end{equation}
for any $\bm{\varpi} \in L^{2}(I;H^{1}(\Omega)^{2}) \cap L^{\infty}(I \times \Omega)^{2}$ and any $n \in \N$, and 
\begin{equation}\label{5-12.55}
\begin{array}{l}
\ds \int_{I} (\alpha_{0}(\overline{\v}_{n}(t))(\widehat{\theta}_{n})_{t}(t), \overline{\theta}_{n}(t) - \zeta(t))_{L^{2}(\Omega)} dt 
    \\[2ex]
    \hspace{7ex} + \Phi_{\nu_{n}, \sigma_{n}}^{I}(\overline{\v}_{n}(t);\overline{\theta}_{n}(t)) \le \Phi_{\nu_{n}, \sigma_{n}}^{I}(\overline{\v}_{n}(t);\zeta(t))
\end{array}
\end{equation}
for any $\zeta \in L^{2}(I;H^{1}(\Omega))$ and any $n \in \N$. 
\medskip

Let us take any $\zeta \in D(\Phi_{\nu_{0}}^{I}(\v;{}\cdot\,)) $. On account of Lemma \ref{lem-gamma}, we can find a sequence $\{ \zeta_{n} \}_{n=1}^{\infty} \subset L^{2}(I ; H^{1}(\Omega))$ such that
\begin{equation*}
\zeta_{n} \to \zeta \mbox{ in } L^{2}(I;L^{2}(\Omega)) \ \ \mbox{ and }\ \ \Phi_{\nu_{n}, \sigma_{n}}^{I}(\overline{\v}_{n};\zeta_{n}) \to \Phi_{\nu_{0}}^{I}(\v;\zeta), 
\end{equation*}
as $n \to \infty$. Then, by (\ref{5-8})-(\ref{5-12}), (\ref{5-12.55}), Lemmas \ref{lem5-1}, \ref{lem-gamma}, and Remark \ref{Rem.Euc}, we see that 
\begin{align}\label{5-13}
    \ds \int_{I} (\alpha_{0} & \ds (\v(t))\theta_{t}(t), \theta(t)-\zeta(t))_{L^{2}(\Omega)}dt + \Phi_{\nu_{0}}^{I}(\v(t);\theta(t)) 
    \nonumber
    \\
    \le & \ds \lim_{n \to \infty} \int_{I} (\alpha_{0}(\overline{\v}_{n}(t))(\widehat{\theta}_{n})_{t}(t), \overline{\theta}_{n}(t) - \zeta(t))_{L^{2}(\Omega)} dt + \liminf_{n \to \infty} \Phi_{\nu_{n}, \sigma_{n}}^{I}(\overline{\v}_{n};\overline{\theta}_{n}) 
    \nonumber
    \\
    \le & \lim_{n \to \infty} \Phi_{\nu_{n}, \sigma_{n}}^{I}(\overline{\v}_{n};\zeta_n) = \Phi_{\nu_{0}}^{I}(\v;\zeta).
\end{align}
Since the choices of the bounded open interval $I \subset \subset (0,\infty)$ and the function $ \zeta \in D(\Phi_{\nu_0}^I(\v;\cdot)) $ is arbitrary, we derive the variational inequality in (S3), as a straightforward consequence of \eqref{5-13}.
\medskip

Next, we put $\zeta = \theta$ in (\ref{5-13}). Then, it can be seen that 
\begin{equation*}
\Phi_{\nu_{0}}^{I}(\v;\theta) \le \liminf_{n \to \infty}\Phi_{\nu_{n}, \sigma_{n}}^{I}(\overline{\v}_{n};\overline{\theta}_{n}) \le \limsup_{n \to \infty}\Phi_{\nu_{n}, \sigma_{n}}^{I}(\overline{\v}_{n};\overline{\theta}_{n}) \le \Phi_{\nu_{0}}^{I}(\v;\theta), 
\end{equation*}
that is 
\begin{equation*}\label{5-14}
\lim_{n \to \infty}\Phi_{\nu_{n}, \sigma_{n}}^{I}(\overline{\v}_{n};\overline{\theta}_{n}) = \Phi_{\nu_{0}}^{I}(\v;\theta).
\end{equation*}
Also, from \eqref{5-9}--\eqref{5-12.1}, (AP2), and Lemma \ref{lem5-0}, one can see that:
    \begin{align*}
        \liminf_{n \to \infty} & \int_I \int_\Omega \alpha(\overline{\v}_n(t)) |\nabla \overline{\theta}_n|_{\sigma_n} dx dt 
        \\
        \geq & \liminf_{n \to \infty} \int_I \int_\Omega \alpha(\overline{\v}_n(t)) \bigl( q_0(\sigma_n) |\nabla \overline{\theta}_n(t)| -r_0(\sigma_n) \bigr) dx dt
        \\
        \geq & \int_I \int_\Omega d \bigl[ \alpha(\v(t)) |D \theta(t)| \bigr] \, dt,
    \end{align*}
    and
    \begin{align*}
        \liminf_{n \to \infty} & \int_I \int_\Omega \beta(\overline{\v}_n(t)) |\nabla (\nu_n \overline{\theta}_n)(t)|^2 dx dt = \liminf_{n \to \infty} \bigl| {\ts \sqrt{\beta(\overline{\v}_n)} |\nabla (\nu_n \overline{\theta}_n)|} \bigr|_{L^{2}(I;L^2(\Omega)^N)}^2
        \\
        \geq & \bigl| {\ts \sqrt{\beta(\v)} |\nabla (\nu_0 \theta)|} \bigr|_{L^{2}(I;L^2(\Omega)^N)}^2 = \int_I \int_\Omega \beta(\v(t)) |\nabla (\nu_0 \theta)(t)|^2 dx dt.
    \end{align*}
The above observations mean that
\begin{equation}\label{5-16}
\int_{I}\int_{\Omega} \alpha(\overline{\v}_{n}(t))|\nabla\overline{\theta}_{n}(t)|_{\sigma_{n}} dxdt \to \int_{I}\int_{\Omega} d[\alpha(\v(t))|D\theta(t)|]dt,
\end{equation}
and 
\begin{equation}\label{5-16.1}
\int_{I}\int_{\Omega} \beta(\overline{\v}_{n}(t))|\nabla(\nu_{n}\overline{\theta}_{n})(t)|^{2} dxdt \to \int_{I}\int_{\Omega} \beta(\v(t))|\nabla(\nu_{0}\theta)(t)|^{2}dxdt
\end{equation}
as $n \to \infty$. By (\ref{5-8}), (\ref{5-11}), (\ref{5-12}), (\ref{5-16}), (A2), and {\bf ($\sharp$b)}, we apply Lemma \ref{lem5-1} with $\rho = \alpha(\v)$, $\{\rho_{n}\}_{n=1}^{\infty} = \{\alpha(\overline{\v}_{n})\}_{n=1}^{\infty}$, $\zeta = \theta$, $\{\zeta_{n}\}_{n=1}^{\infty} = \{\overline{\theta}_{n}\}_{n=1}^{\infty}$, $\omega =1$, and $\{\omega_{n}\}_{n=1}^{\infty} = \{1\}$. Then, we have 
\begin{equation}\label{5-16.2}
\int_{I}\int_{\Omega} |\nabla \overline{\theta}_{n}(t)|_{\sigma_{n}}dxdt \to \int_{I}\int_{\Omega} |D \theta(t)| dt 
\end{equation}
as $n \to \infty$. 
\medskip

Besides, we take $I = (t_{0},t_{1}) \subset \R$ with $ 0 < t_0 < t_1 < \infty $. Using (A2), (\ref{5-2}), (\ref{5-5}), (\ref{5-7.5}), and {\bf ($\sharp$b)}, it is deduced that 
\begin{subequations}\label{5-16.5}
\begin{align}\label{5-16.5a}
    & |\overline{\theta}_n -\underline{\theta}_n|_{L^\infty(I; L^2(\Omega))} \leq {\ts \sqrt{h_n}} |(\widehat{\theta}_n)_t|_{L^2(I; L^2(\Omega))} \to 0,
\end{align}
\begin{align}\label{5-16.5b}
    \left| \rule{-1pt}{16pt} \right. \int_I \int_\Omega & |\nabla \overline{\theta}_n(t)|_{\sigma_n} dx dt -\int_I \int_\Omega |\nabla \underline{\theta}_n(t)|_{\sigma_n} dxdt \left. \rule{-1pt}{16pt} \right|
    \nonumber
    \\
    \leq &~ h_n \left| \int_\Omega \bigl( |\nabla \overline{\theta}_n(t_1)|_{\sigma_n} -|\nabla \underline{\theta}_n(t_0)|_{\sigma_n} \bigr) dx \right|
    \nonumber
    \\
    \leq &~ \frac{2 h_n}{\delta_*} \left( \F_{\nu_n}(\v_{0, n}, \theta_{0, n}) +c |u^\dag|_{L^1(\Omega)} +c^2 |\overline{u}_n -u^\dag|_{L^2(I; L^2(\Omega))}^2 \right) \to 0,
\end{align}
\end{subequations}
as $n \to \infty$. As is seen from the above convergences, (\ref{5-8}), (\ref{5-11}), (\ref{5-12}), (A2), and {\bf ($\sharp$b)}, we take any $\bm{\varpi} \in [H^{1}(\Omega) \cap L^{\infty}(\Omega)]^{2}$, and apply Lemma \ref{lem5-1} with $\rho = 1$, $\{\rho_{n}\}_{n=1}^{\infty} = \{ 1 \}$, $\zeta = \theta$, $\{\zeta_{n}\}_{n=1}^{\infty} = \{ \underline{\theta}_{n}\}_{n=1}^{\infty}$, $\omega= \bm{\varpi} \cdot [\nabla\alpha](\v)$, and $\{ \omega_{n} \}_{n=1}^{\infty} = \{ \bm{\varpi} \cdot [\nabla \alpha](\overline{\v}_{n})\}$. Then, we see that
\begin{equation}\label{5-17}
\lim_{n \to \infty} \int_{I}\int_{\Omega} \bm{\varpi} \cdot [\nabla\alpha](\overline{\v}_{n}(t))|\nabla \underline{\theta}_{n}(t)|_{\sigma_{n}} dxdt = \int_{I}\int_{\Omega} d[\bm{\varpi} \cdot [\nabla\alpha](\v(t))|D\theta(t)| ]dt
\end{equation}
for any $\bm{\varpi} \in [H^{1}(\Omega) \cap L^{\infty}(\Omega)]^{2}$. 

On the other hand, (\ref{5-12.1}), (\ref{5-16.1}), and the uniform convexity of $L^{2}$-based topology imply that
\begin{equation*}
\ts\sqrt{\beta(\overline{\v}_{n})} \nabla (\nu_{n}\overline{\theta}_{n}) \to \sqrt{\beta(\overline{\v})} \nabla (\nu_{0}\theta)\ \ \mbox{ in } L^{2}(I;L^{2}(\Omega)^{N}), 
\end{equation*}
and hence 
\begin{equation}\label{5-17.5}
\nabla(\nu_{n}\overline{\theta}_{n}) \to \nabla(\nu_{0}\theta)\ \ \mbox{ in } L^{2}(I;L^{2}(\Omega)^{N}), 
\end{equation}
as $n \to \infty$. In addition,  by (A2), (\ref{5-5}), (\ref{5-7.4}), (\ref{5-7.5}), and {\bf ($\sharp$b)}, we can show that 
\begin{align}\label{5-17.6}
    \ds \left| \rule{-1pt}{16pt} \right. \int_{I}\int_{\Omega} & |\nabla (\nu_{n}\overline{\theta}_{n})|^{2} dxdt - \int_{I}\int_{\Omega} |\nabla (\nu_{n}\underline{\theta}_{n})|^{2} dxdt \left. \rule{-2pt}{16pt} \right| 
    \nonumber
    \\
    \leq &~ \frac{2 h_n}{\delta_*} \left( \F_{\nu_n}(\v_{0, n}, \theta_{0, n}) +c |u^\dag|_{L^1(\Omega)} +c^2 |\overline{u}_n -u^\dag|_{L^2(I; L^2(\Omega))}^2 \right) \to 0,
\end{align}
as $n \to \infty$. As a consequence of (\ref{5-8}), (\ref{5-9}), (\ref{5-12.1}), (\ref{5-17.5}), (\ref{5-17.6}), and {\bf ($\sharp$b)}, it is observed that
\begin{equation}\label{5-17.7}
\left\{ \begin{array}{l}
\ds \nu_{n}\overline{\theta}_{n} \to \nu_{0}\theta,\ \ \nu_{n}\underline{\theta}_{n} \to \nu_{0}\theta \ \ \mbox{ in } L^{2}(I; H^{1}(\Omega)), \vspace{3mm}\\
\ds \sqrt{{\bm \varpi} \cdot [\nabla\beta](\overline{\v}_{n})}\nabla(\nu_{n}\underline{\theta}_{n}) \to \sqrt{{\bm \varpi} \cdot [\nabla \beta](\v)}\nabla(\nu_{0}\theta)\ \mbox{ in } L^{2}(I;L^{2}(\Omega)^{N})
\end{array} \right.
\end{equation}
for any ${\bm \varpi} \in [H^{1}(\Omega) \cap L^{\infty}(\Omega)]^{2}$ as $n \to \infty$. 

From (\ref{5-8})--(\ref{5-12}), (\ref{5-17}), and (\ref{5-17.7}), letting $n \to \infty$ in (\ref{5-12.5}) implies 
\begin{equation*}
\begin{array}{l}
\ds \int_{I} (\v_{t}(t), \v(t) - \bm{\varpi})_{L^{2}(\Omega)^{2}} dt + \int_{I} (\nabla \v(t), \nabla (\v(t) - \bm{\varpi}))_{L^{2}(\Omega)^{2N}} dt \vspace{3mm}\\
\ds \hspace{5mm} + \int_{I} ([\nabla G](u;\v)(t), \v(t) - \bm{\varpi})_{L^{2}(\Omega)^{2}} dt + \int_{I}\int_{\Omega} \gamma (\v(t)) dxdt \vspace{3mm}\\
\ds \hspace{5mm} + \int_{I} \int_{\Omega} d[(\v(t) - \bm{\varpi}) \cdot [\nabla\alpha](\v(t)) |D\theta(t)|] dt \vspace{3mm}\\
    \ds \hspace{5mm} + \int_{I} \int_{\Omega} (\v(t) - \bm{\varpi}) \cdot [\nabla\beta](\v(t)) |\nabla (\nu_{0} \theta)(t)|^{2} dxdt \ds \le \int_{I}\int_{\Omega} \gamma(\bm{\varpi}) dxdt
\end{array}
\end{equation*}
for any $\bm{\varpi} \in [H^{1}(\Omega) \cap L^{\infty}(\Omega)]^{2}$, and any $\nu_{0} \in [0,1)$. Since the open interval $I \subset (0,\infty)$ is arbitrary, $\v = [w,\eta]$ satisfies the variational inequalities in (S2). 

We next check the initial condition. To see this, we fix $t \in (0,\infty)$, $\ell, n \in \N$. We take $[h,\nu,\sigma] = [h_{n},\nu_{n},\sigma_{n}]$ and $[\bm{w}_{0},\omega_{0}] = [\v_{0,\ell},\theta_{0,\ell}]$ in (\ref{5-7}). By (A2), (A3), (A4), 
(\ref{enrgyF_nu}), (\ref{5-1}), (\ref{5-7.4})-(\ref{5-12.1}), and Lemma \ref{lem-gamma}, letting $n \to \infty$ gives that 
\begin{equation*}
\begin{array}{l}
\ds \frac{1}{2} (|\v(t) - \v_{0,\ell}|_{L^{2}(\Omega)^{2}}^{2} + A_{\ast}|\theta(t) - \theta_{0,\ell}|_{L^{2}(\Omega)}^{2}) \vspace{3mm}\\
\ds \hspace{5mm} + \frac{B_{\ast}}{4} \int_{0}^{t} |\nabla \v|_{L^{2}(\Omega)^{2N}}^{2} dt + \frac{B_{\ast}}{2} \delta_{\ast
    } ||D\theta(\cdot)|(\Omega)|_{L^{1}(0,t)} + \frac{B_{\ast}}{2} \delta_{\ast} |\nabla(\nu_0 \theta)|_{L^{2}(0,t;L^{2}(\Omega)^{N})}\vspace{3mm}\\
\ds \le \frac{1}{2} (|\v_{0} - \v_{0,\ell}|_{L^{2}(\Omega)^{2}}^{2} + A_{\ast}|\theta_{0} - \theta_{0,\ell}|_{L^{2}(\Omega)}^{2}) + 2tC_{\ast} (1 + |\v_{0,\ell}|_{H^{1}(\Omega)^{2}}^{2} + |\theta_{0,\ell}|_{H^{1}(\Omega)}^{2}) \vspace{3mm}\\
\ds \hspace{5mm} + \frac{c^{2}}{2}\int_{0}^{t}|u|_{L^{2}(\Omega)}^{2} d\tau.
\end{array}
\end{equation*}
Hence, it holds that $\v \in L^{2}_\mathrm{loc}([0,\infty);H^{1}(\Omega)^{2})$, $|D\theta(\cdot)|(\Omega) \in L^{1}_\mathrm{loc}([0,\infty))$ and $\nu_0 \theta \in L^{2}_\mathrm{loc}([0,\infty);H^{1}(\Omega))$ hold. 

Furthermore, it follows that 
\begin{equation*}
\begin{array}{l}
\ds \frac{1}{2} (|\v(t) - \v_{0}|_{L^{2}(\Omega)^{2}}^{2} + A_{\ast}|\theta(t) - \theta_{0}|_{L^{2}(\Omega)}^{2}) \vspace{3mm}\\
\ds \le |\v(t) - \v_{0,\ell}|_{L^{2}(\Omega)^{2}}^{2} + A_{\ast}|\theta(t) - \theta_{0,\ell}|_{L^{2}(\Omega)}^{2} + |\v_{0,\ell} - \v_{0}|_{L^{2}(\Omega)}^{2} + A_{\ast}|\theta_{0,\ell} - \theta_{0}|_{L^{2}(\Omega)}^{2}, 
\end{array}
\end{equation*}
for any $t \in (0,\infty)$ and $\ell \in \N$. Combining the above two inequalities, we can deduce that 
\begin{equation*}
\begin{array}{l}
\ds \limsup_{t \downarrow 0} (|\v(t) - \v_{0}|_{L^{2}(\Omega)^{2}}^{2} + A_{\ast}|\theta(t) - \theta_{0}|_{L^{2}(\Omega)}^{2}) \vspace{3mm}\\
\ds \le 4(|\v_{0,\ell} - \v_{0}|_{L^{2}(\Omega)^{2}}^{2} + A_{\ast}|\theta_{0,\ell} - \theta_{0}|_{L^{2}(\Omega)}^{2}).
\end{array}
\end{equation*}
By (\ref{5-1}), (\ref{5-7.51}), and (\ref{5-11})--(\ref{5-12}), the above inequality implies $[\v,\theta] \in C([0,\infty);L^{2}(\Omega)^{3})$ and $[\v(0),\theta(0)] = [\v_{0},\theta_{0}]$ in $L^{2}(\Omega)^{3}$. Therefore, (S1) is verified.
\medskip

Finally, we prove the energy dissipation (S4). The verification  of (S4) is reduced to the following key lemma.
\begin{lem}[Energy inequality]\label{key-lem}
Let $\nu_{\ast}$ be the positive constant, obtained in Lemma \ref{lem4-1}, and let $\nu_{0} \in [0,\nu_{\ast})$ be a fixed constant. Let  $[\v(t), \theta(t)] = [w(t), \eta(t), \theta(t)]$ be the triplet of functions, as in (\ref{5-8})--(\ref{5-12.1}). Let $\J^{{u}^\dag}_{\ast} \in BV_\mathrm{loc}((0,\infty))$ be the function, as in \eqref{5-12.3}. Then, 
\begin{align}\label{5-18}
    \J^{{u}^\dag}_{\ast}(t) &~ = \widehat{\F}_{\nu_{0}}^{{u}^\dag}(t):= \F_{\nu_{0}}(\v(t),\theta(t)) 
    \nonumber
    \\
    &~ + c \int_{\Omega} {u}^\dag w(t) dx - c^{2} \int_{0}^{t} \int_{\Omega} |u(\tau) - {u}^\dag|^{2} dxd\tau, 
\end{align}
for a.e. $t \in (0,\infty)$. Moreover, the following energy inequality holds:
\begin{equation}\label{5-30}
\frac{1}{2} \int_{s}^{t} |\v_{t}(\tau)|_{L^{2}(\Omega)^{2}}^{2} d\tau + \int_{s}^{t} |\sqrt{\alpha_{0}(\v(\tau))} \theta_{t}(\tau)|_{L^{2}(\Omega)}^{2} d\tau + \J^{{u}^\dag}_{\ast}(t) \le \J^{{u}^\dag}_{\ast}(s)
\end{equation}
for a.e. $0 < s < t < \infty$. Therefore, $\J^{{u}^\dag}_{\ast}$ is non-increasing on $(0,\infty)$. 
\end{lem}
\textbf{Proof.} 
We take any bounded open interval $I \subset \subset (0,\infty)$ and a sequence $\{\v_{n}\}_{n=1}^{\infty} \subset C^{\infty}(\overline{I \times \Omega})^{2}$ such that $\v_{n} \to \v$ in $L^{2}(I;H^{1}(\Omega)^{2})$ as $n \to \infty$. Moreover, we choose $\bm{\varpi} = \v_{n}$ in (\ref{5-12.5}). By (\ref{5-16.2})-(\ref{5-17}), Lemma \ref{lem5-1} with $\rho = 1$, $\{\rho_{n}\}_{n=1}^{\infty} = \{1\}$, $\zeta = \theta$, $\{\zeta_{n}\}_{n=1}^{\infty} = \{ \underline{\theta}_{n}\}_{n=1}^{\infty}$, $\omega = 0$, and $\{\omega_{n}\}_{n=1}^{\infty} = \{ (\overline{\v}_{n} - \v_{n}) \cdot [\nabla\alpha](\overline{\v}_{n})\}_{n=1}^{\infty}$, and (\ref{5-17.7}) with ${\bm \varpi} = \overline{\v}_{n}(t) - \v_{n}(t)$ follow that 
\begin{equation}\label{5-19}
\begin{array}{l}
\ds \int_{I} |\nabla \v(t)|_{L^{2}(\Omega)^{2N}}^{2} dt \le \liminf_{n \to \infty} \int_{I} |\nabla \overline{\v}_{n}(t)|_{L^{2}(\Omega)^{2N}}^{2} dt \le \limsup_{n \to \infty} \int_{I} |\nabla \overline{\v}_{n}(t)|_{L^{2}(\Omega)^{2N}}^{2} dt \vspace{3mm}\\
    \ds \le \lim_{n \to \infty} \left\{ \int_{I}|\nabla \v_{n}(t)|_{L^{2}(\Omega)^{2N}}^{2} dt - 2 \int_{I} ( (\widehat{\v}_{n})_{t}(t) + [\nabla G](\overline{u}_{n};\overline{\v}_{n})(t), (\overline{\v}_{n} - \v_{n})(t) )_{L^{2}(\Omega)^{2}} dt \right. \vspace{3mm}\\
\ds \hspace{5mm} - 2 \int_{I}\int_{\Omega} (\overline{\v}_{n}(t) - \v_{n}(t)) \cdot ([\nabla\alpha](\overline{\v}_{n}(t))|\nabla\underline{\theta}_{n}|_{\sigma_{n}} + \nu_{n}^{2}[\nabla\beta](\overline{\v}_{n}(t))|\nabla\underline{\theta}_{n}(t)|^{2}) dxdt \vspace{3mm}\\
\ds \hspace{5mm} \left. + 2 (\int_{I}\int_{\Omega}\gamma(\v_{n}(t))dxdt - \int_{I}\int_{\Omega}\gamma(\overline{\v}_{n}(t))dxdt) \right\} 
\ds = \int_{I} |\nabla \v(t)|_{L^{2}(\Omega)^{2N}}^{2} dt. 
\end{array}
\end{equation}
Using (\ref{5-11})-(\ref{5-12.1}), (\ref{5-16}), (\ref{5-17.7}), (\ref{5-19}), \textbf{($\sharp$e)}, and the uniform convexities of the $L^{2}$-type topologies, we see that 
\begin{equation}\label{5-19.5}
\left\{ \begin{array}{l}
\ds \overline{\v}_{n} \to \v \mbox{ in } L^{2}(I;H^{1}(\Omega)^{2}), 
    \\[1ex]
\ds \int_{I} \overline{\F}_{n}^{{u}^\dag}(t) dt \to \int_{I} \widehat{\F}^{{u}^\dag}_{\nu_{0}}(t) dt, 
\end{array} \right. 
\end{equation}
as $n \to \infty$. Here, putting:
\begin{equation*}
    \begin{cases}
        C_0 := c|u^\dag|_{L^1(\Omega)} +c^2 |u -u^\dag|_{L^2(I; L^2(\Omega))}^2,
        \\[1ex]
        C_n := c|u^\dag|_{L^1(\Omega)} +c^2 |\overline{u}_n -u^\dag|_{L^2(I; L^2(\Omega))}^2 \mbox{ for } n = 1, 2, 3, \dots,
    \end{cases}
\end{equation*}
we easily see from \eqref{5-8}, \eqref{5-7.52}, \eqref{5-18}, and \textbf{($\sharp$a)} that:
\begin{equation}\label{ken5-18-01}
    \begin{cases}
        \widehat{\F}_{\nu_0}^{u^\dag} +C_0 \geq 0, \mbox{ and } \underline{\F}_{n}^{u^\dag} +C_n \geq 0 \mbox{ for $ n = 1, 2, 3, \dots $,}
        \\[1ex]
        C_n \to C_0 \mbox{ as $ n \to \infty $.}
    \end{cases}
\end{equation}
\eqref{5-7.4}--\eqref{5-7.5}, \eqref{5-11}--\eqref{5-12.1}, \eqref{ken5-18-01}, and \textbf{($\sharp$e)} enable us to compute:
\begin{align}\label{5-20}
    \left| \rule{-1pt}{16pt} \right. \int_I \bigl( \underline{\F}_n^{u^\dag} & (t) + C_n \bigr) dt -\int_I \bigl( \widehat{\F}_{\nu_0}^{u^\dag}(t) + C_0 \bigr) dt \left. \rule{-2pt}{16pt} \right|
    \nonumber
    \\
    \leq &~ \left| \int_I \underline{\F}_n^{u^\dag}(t) dt -\int_I \overline{\F}_n^{u^\dag}(t) dt \right| 
    \nonumber
    \\
    & \quad +\left| \int_I \bigl( \overline{\F}_n^{u^\dag}(t) +C_n \bigr) dt -\int_I \bigl( \widehat{\F}_{\nu_0}^{u^\dag}(t) +C_0 \bigr) dt \right|
    \nonumber
    \\
    \leq &~ 2 h_n \bigl( \F_{\nu_n, \sigma_n}(\v_{0, n}, \theta_{0, n}) +c|u^\dag|_{L^1(\Omega)} +c^2 |\overline{u}_n -u^\dag|_{L^2(I; L^2(\Omega))}^2 \bigr)
    \nonumber
    \\
    & \quad +\left|  \int_I \overline{\F}_n^{u^\dag}(t) dt -\int_I \widehat{\F}_{\nu_0}^{u^\dag}(t) dt \right| +|C_n -C_0| \mathscr{L}^{1}(I)
    \nonumber
    \\
    \to &~ 0,
\end{align}
as $ n \to \infty $, for any bounded open interval $ I \subset \subset (0, \infty) $.

For any given bounded open set $A \subset (0,\infty)$, we denote by $\mathcal{I}_{A}$ the at most countable class of pointwise disjoint open intervals such that $\cup_{\tilde{I} \in \mathcal{I}_{A}} \tilde{I} = A$. Here (\ref{5-20}) yields that 
\begin{equation*}
    \sum_{\tilde{I} \in \tilde{\mathcal{I}}} \int_{\tilde{I}} \bigl( \widehat{\F}^{{u}^\dag}_{\nu_{0}}(t) +C_0 \bigr) dt \le \liminf_{n \to \infty} \int_{A} \bigl( \underline{\F}_{n}^{{u}^\dag} (t) +C_n \bigr) dt,
\end{equation*}
for any finite subclass $\tilde{\mathcal{I}} \subset \mathcal{I}_{A}$. Here, it follows that
\begin{equation}\label{5-21}
\int_{A} \bigl( \widehat{\F}^{{u}^\dag}_{\nu_{0}}(t) +C_0 \bigr) dt \le \liminf_{n \to \infty} \int_{A} \bigl( \underline{\F}_{n}^{{u}^\dag} (t) +C_n \bigr) dt,
\end{equation}
for all $A \subset \subset (0,\infty)$. Applying \cite[Proposition 1.80]{AFP} to (\ref{5-19.5})-(\ref{5-21}), we can see that 
\begin{equation}\label{5-17'}
    \begin{array}{c}
        \underline{\F}_{n}^{{u}^\dag} = \bigl( \underline{\F}_{n}^{{u}^\dag} +C_n \bigr) -C_n \to \bigl( \widehat{\F}^{{u}^\dag}_{\nu_{0}}+C_0 \bigr) -C_0 = \widehat{\F}^{{u}^\dag}_{\nu_{0}},
        \\[1ex]
        \mbox{ weakly-$*$ in } \mathcal{M}_\mathrm{loc}((0,\infty)),
    \end{array}
\end{equation}
as $n \to \infty$. Consequently, it follows from \eqref{5-12.3} and \eqref{5-17'} that
\begin{equation}\label{5-29}
    \widehat{\F}^{{u}^\dag}_{\nu_{0}}(t) = \J^{{u}^\dag}_{\ast}(t)
\end{equation}
a.e. $t \in (0,\infty)$. 

In addition, employing \eqref{5-12.3} and \eqref{5-29}, and passing to the limit as $n \to \infty$ in (\ref{5-5}) with $[h,\nu,\sigma] = [h_{n},\nu_{n},\sigma_{n}]$, the energy inequality (\ref{5-30}) holds. Also, \eqref{5-8}, \eqref{5-18}, and \eqref{5-29} show that 
\begin{equation}\label{5-31}
    \J_*^{u^\dag} = \widehat{\F}^{{u}^\dag}_{\nu_{0}} \in L^{1}_\mathrm{loc}([0,\infty)) \cap L^{\infty}_\mathrm{loc}((0,\infty)). 
\end{equation}
Combining (\ref{5-30}), (\ref{5-29}), and (\ref{5-31}), the desired assertion holds. 
\hfill $\Box$ 
\begin{rem}\label{left-conti}
\begin{em}
In the following, we assign the left-continuous expression of $t \in (0,\infty) \mapsto \hat{\F}_{\nu_{0}}^{{u}^\dag}(\v(t),\theta(t))$ to the function $\J^{{u}^\dag}_{\ast}(t)$ in (\ref{5-18}). Then, (\ref{5-30}) can be satisfied for all $0 < s \le t < \infty$ by the nonincreasing property of $\J^{{u}^\dag}_{\ast}$. 
\end{em}
\end{rem}

\section{Proof of Main Theorem 2}
Let $\nu_{0} \ge 0$ be a fixed constant. In this section, we prove the large-time behavior for the solutions to (S)$_{\nu_{0}}$ which are constructed in Section 5. To prove this, we recall the results for weighted total variations.
\begin{lem}[cf. {\cite[Lemma 4.4]{MSW}}]\label{lem6-1}
Let $I \subset (0,\infty)$ be a bounded open interval. Let $\rho \in C(\overline{I}; L^{2}(\Omega)) \cap L^{\infty}(I; H^{1}(\Omega)) \cap L^{\infty}(I \times \Omega)$ and $\{\rho_{n}\}_{n=1}^{\infty} \subset C(\overline{I}; L^{2}(\Omega)) \cap L^{\infty}(I; H^{1}(\Omega))$, $\zeta \in L^{2}(I;L^{2}(\Omega))$, and $\{\zeta_{n}\}_{n=1}^{\infty} \subset L^{2}(I;L^{2}(\Omega))$ be such that $|D\zeta(\cdot)|(\Omega) \in L^{1}(I)$, and 
\begin{equation*}
\left\{ 
\begin{array}{l}
\ds \rho(t),\, \rho_{n}(t) \in W_{c}(\Omega) \mbox{ a.e. } t \in I, \vspace{3mm}\\
\ds t \in I \mapsto \int_{\Omega} d[\rho(t)|D\zeta(t)|] \mbox{ and } t \in I \mapsto \int_{\Omega} d[\rho_{n}(t)|D\zeta_{n}(t)|], \mbox{ for } n \in \N, \mbox{ are measurable,} \vspace{3mm}\\
\ds \rho_{n}(t) \to \rho(t) \mbox{ in } L^{2}(\Omega) \mbox{ and weakly in } H^{1}(\Omega) \mbox{ a.e. } t \in I, \mbox{ as } n \to \infty, \vspace{3mm}\\
\ds \zeta_{n}(t) \to \zeta(t) \mbox{ a.e. } t \in I, \mbox{ as } n \to \infty, \vspace{3mm}\\
\ds \rho \ge \delta_{0} \mbox{ and } \inf_{n \in \N} \rho_{n} \ge \delta_{0} \mbox{ a.e. in } I \times \Omega, \mbox{ for some constant } \delta_{0}>0 
\end{array}
\right. 
\end{equation*}
are satisfied. Also, let $\omega \in C(\overline{I}; L^{2}(\Omega)) \cap L^{\infty}(I; H^{1}(\Omega)) \cap L^{\infty}(I \times \Omega)$ and $\{\rho_{n}\}_{n=1}^{\infty} \subset C(\overline{I}; L^{2}(\Omega)) \cap L^{\infty}(I; H^{1}(\Omega)) \cap L^{\infty}(I \times \Omega)$ be such that 
\begin{equation*}
\left\{ 
\begin{array}{l}
\ds \omega_{n}(t) \to \omega(t) \mbox{ in } L^{2}(\Omega) \mbox{ and weakly in } H^{1}(\Omega), \mbox{ a.e. } t\in I, \mbox{ as } n \to \infty, \vspace{3mm}\\
\ds |\omega| \le M_{0} \mbox{ and } \sup_{n \in \N} |\omega_{n}| \le M_{0} \mbox{ a.e. in } I\times\Omega, \mbox{ for some constant } M_{0} > 0
\end{array}
\right. 
\end{equation*}
hold. In addition, let us assume 
\begin{equation*}
\int_{I}\int_{\Omega} d[\rho_{n}(t)|D\zeta_{n}(t)|] dt \to \int_{I} \int_{\Omega} d[\rho(t)|D\zeta(t)|]dt
\end{equation*}
as $n\to\infty$. Then, 
\begin{equation*}
\int_{I}\int_{\Omega} d[\omega_{n}(t)|D\zeta_{n}(t)|] dt \to \int_{I} \int_{\Omega} d[\omega(t)|D\zeta(t)|]dt
\end{equation*}
as $n\to\infty$.
\end{lem}
\begin{lem}[cf. {\cite[Theorem 4.1]{SW15}}]\label{lem-gamma2}
Let $I \subset (0,\infty)$ be any open interval. If a function $ \tilde{\v} \in C(\overline{I}; L^2(\Omega)^2) \cap L^\infty(I; H^1(\Omega)^2) \cap L^\infty(I \times \Omega)^2 $ and a sequence $ \{ \tilde{\v}_n\}_{n = 1}^\infty \subset C(\overline{I}; L^2(\Omega)^2) \cap L^\infty(I; H^1(\Omega)^2) \cap L^\infty(I \times \Omega)^2 $ satisfy that $ \{ \tilde{\v}_n \}_{n = 1}^\infty $ is bounded in $ L^\infty(I \times \Omega)^{2} $, and $ \tilde{\v}_n(t) \to \tilde{\v}(t) $ in $ L^2(\Omega)^2 $ and weakly in $ H^1(\Omega)^2 $ as $ n \to \infty $, for a.e. $ t \in I $, then the sequence $ \{ \Phi_{\nu_{0}}^I(\tilde{\v}_n;{}\cdot{}) \}_{n = 1}^\infty $ converges to $ \Phi_{\nu_{0}}^I(\tilde{\v};{}\cdot{}) $ on $ L^2(I; L^2(\Omega)) $, in the sense of $ \Gamma $-convergence, as $ n \to \infty $.
\end{lem}
\begin{rem}[cf. {\cite[Remark 4.1]{MSW}} ]\label{rem6-1}
\begin{em}
    If $\rho \in W_{0}(\Omega)$, $\{\rho_{n}\}_{n=1}^{\infty} \subset W_{0}(\Omega)$, $\zeta \in L^{2}(\Omega)$, and $\{\zeta_{n}\}_{n=1}^{\infty} \subset L^{2}(\Omega)$ fulfill that 
\begin{equation*}
\begin{array}{l}
\ds \rho_{n} \to \rho \mbox{ in } L^{2}(\Omega) \mbox{ and weakly in } H^{1}(\Omega),\hspace{3mm} 
\ds \zeta_{n} \to \zeta \ \mbox{ in } L^{2}(\Omega)
\end{array}
\end{equation*}
as $n \to \infty$, and 
\begin{equation*}
\rho \in W_{c}(\Omega) \mbox{ or } \sup_{n \in \N} |D\zeta_{n}|(\Omega) < \infty, 
\end{equation*}
then, 
\begin{equation*}
\liminf_{n \to \infty} \int_{\Omega} d[\rho_{n}|D\zeta_{n}|] \ge \int_{\Omega} d[\rho|D\zeta|]. 
\end{equation*}
\end{em}
\end{rem}

We fix the function $u_{\infty}$ which is defined in (A6), and set ${u}^\dag = u_{\infty}$ in Lemma \ref{key-lem}. Then, the functional $\J_{\ast}^{u_{\infty}}$ is nonincreasing on $(0, \infty)$, and satisfies (\ref{5-18}) and (\ref{5-30}) with ${u}^\dag = u_{\infty}$. In view of (S1), Remarks \ref{left-conti} and \ref{rem6-1}, we see that
\begin{equation}\label{6-1}
\begin{array}{l}
\ds \frac{1}{2} \int_{s}^{t} |\v_{t}(\tau)|_{L^{2}(\Omega)^{2}}^{2} d\tau + \int_{s}^{t} |\sqrt{\alpha_{0}(\v(\tau))} \theta_{t}(\tau)|_{L^{2}(\Omega)}^{2} d\tau + \widehat{\F}^{u_{\infty}}_{\nu_{0}}(t) \le \J^{u_{\infty}}_{\ast}(s), 
\end{array}
\end{equation}
for all $ 1 \leq  s \le t < \infty$.
\medskip

Let $[\v,\theta] \in C([0,\infty);L^{2}(\Omega)^{3}) \cap W^{1,2}_\mathrm{loc}((0,\infty);L^{2}(\Omega)^{3})$ be an energy dissipative solution to (S)$_{\nu_{0}}$. 
Then, from (\ref{6-1}), (A2), and Remark \ref{Rem.Note05} (Fact\,4), we estimate that:
\begin{align*}
    \frac{1}{2} \int_1^t & \bigl( |\v_t(\tau)|_{L^2(\Omega)}^2 +\delta_* |\theta_t(\tau)|_{L^2(\Omega)}^2 \bigr) dt +\frac{1}{2} |\nabla \v(t)|_{L^2(\Omega)^N}^2 +\delta_* |D \theta(t)|(\Omega)
    \\
    \leq & \F_{\nu_0}(\v(1 -0), \theta(1 -0)) +c|u_\infty|_{L^1(\Omega)} +c^2 |u -u_\infty|_{L^2(1, \infty; L^2(\Omega))}^2 =: K_*,
\end{align*}
for all $ t \geq 1 $. This 
implies that $[\v_{t},\theta_{t}] \in L^{2}(1,\infty;L^{2}(\Omega)^{3})$, and hence
\begin{equation}\label{6-1.4}
    \ds [\v_{t}(\cdot+s), \theta_{t}(\cdot+s)] \to 0 ~ (= [0, 0, 0]) \mbox{ in } L^{2}((0,1);L^{2}(\Omega)^{3})
\end{equation}
as $s \to \infty$, and 
\begin{equation}\label{6-1.6}
    \ds \{[\v(t),\theta(t)]\ |\ t \ge 1\} \subset F_{1} := \left\{ \begin{array}{l|l} 
        [\tilde{\v}, \tilde{\theta}] \in D_0 & \parbox{5.5cm}{
            $ |\tilde{\theta}| \leq |\theta_0|_{L^\infty(\Omega)} $ a.e. in $ \Omega $, and $ |\nabla \tilde{\v}|_{L^2(\Omega)^N}^2 +\delta_* |D \tilde{\theta}|(\Omega) \leq 2 K_* $
        }
    \end{array} \right\}. 
\end{equation}
Due to (\ref{free-energy}), (A1)--(A6), and Remark \ref{rem6-1}, $F_{1}$ is closed in $L^{2}(\Omega)^{3}$, and bounded in $[H^{1}(\Omega) \cap L^{\infty}(\Omega)]^{2} \times [BV(\Omega) \cap L^{\infty}(\Omega)]$. Consequently, $F_{1}$ is compact in $L^{2}(\Omega)^{3}$. 

From the above, there exist a triplet $[\v_{\infty},\theta_{\infty}] \in L^{2}(\Omega)^{3}$ and a sequence $\{t_{n}\}$ with $1 \le t_{1} < t_{2} < \cdots < t_{n} \to \infty$ as $n \to \infty$ such that 
\begin{equation}\label{6-2}
[\v(t_{n}), \theta(t_{n})] \to [\v_{\infty}, \theta_{\infty}] \mbox{ in } L^{2}(\Omega)^{3}
\end{equation}
as $n \to \infty$. This means that $\omega(\v, \theta) \neq \emptyset$. Moreover, the compactness of $\omega(\v,\theta)$ is given by the compactness of $F_{1}$, and 
\begin{equation*}\label{6-3}
\omega(\v,\theta) = \bigcap_{s \ge 0} \overline{\{[\v(t),\theta(t)]\ |\ t \ge s \}} \subset \overline{\{[\v(t),\theta(t)]\ |\ t \ge 1 \}} \subset F_{1}. 
\end{equation*}

Next, we take any $[\v_{\infty},\theta_{\infty}] \in \omega(\v,\theta)$, with the divergent sequence $ \{ t_n \}_{n = 1}^\infty \subset (0, \infty) $ as in \eqref{6-2}.Then, (i-a) is verified as a direct consequence of (S1), (\ref{6-1.6}), and (\ref{6-2}). 
Also, (\ref{6-1.4})-(\ref{6-1.6}) ensure that 
\begin{equation}\label{6-4}
    \left\{ \hspace{-2ex}
    \parbox{12cm}{\vspace{-2ex}
        \begin{itemize}
            \item $ \{\v_{n}\}_{n=1}^{\infty} := \{\v(\cdot + t_{n})\}_{n=1}^{\infty} $ is  bounded in  $ W^{1,2}(0,1;L^{2}(\Omega)^{2}) \cap L^{\infty}(0,1;H^{1}(\Omega)^{2}) $,
            \item $ \{\theta_{n}\}_{n=1}^{\infty} := \{\theta(\cdot + t_{n})\}_{n=1}^{\infty} $  is  bounded in $ W^{1,2}(0,1;L^{2}(\Omega)) $, and $ \{|D\theta_{n}(\cdot)|(\Omega)\}_{n=1}^{\infty} $ is  bounded in $ L^{\infty}(0,1) $,
            \item $ \{ [\v_{n}(t),\theta_{n}(t)] \ |\ t \in [0,1],\ n \in \N\} \subset D_{\nu_{0}}^{\ast}(\theta_{0}) $. 
        \vspace{-2ex}
        \end{itemize}
    }
    \right.
\end{equation}
From (\ref{6-1.4}), (\ref{6-4}), and the compactness results as in \cite[Chapter 3]{AFP} and \cite[Corollary 4]{Simon}, we can obtain that 
\begin{equation}\label{6-5}
\left\{ 
\begin{array}{l}
\ds (\v_{n})_{t} \to 0 \mbox{ in } L^{2}(0,1;L^{2}(\Omega)^{2}), \mbox{ and } (\theta_{n})_{t} \to 0 \mbox{ in } L^{2}(0,1;L^{2}(\Omega)), \vspace{3mm}\\
\ds \v_{n} \to \v \mbox{ in } W^{1,2}(0,1;L^{2}(\Omega)^{2}), \mbox{ weakly-}\ast \mbox{ in } L^{\infty}(0,1;H^{1}(\Omega)^{2}), \vspace{3mm}\\
\ds \ \ \mbox{ and weakly-}\ast \mbox{ in } L^{\infty}((0,1)\times\Omega)^{2}, \vspace{3mm}\\
\ds \theta_{n} \to \theta_{\infty} \mbox{ in } W^{1,2}(0,1;L^{2}(\Omega)), \mbox{ and weakly-}\ast \mbox{ in } L^{\infty}((0,1)\times\Omega), \vspace{3mm}\\ 
\ds \theta_{n}(t) \to \theta_{\infty} \mbox{ weakly-}\ast \mbox{ in } BV(\Omega), \mbox{ for any } t \in (0,1), 
\end{array}
\right.
\end{equation}
as $n \to \infty$, by taking subsequences (not relabeled) if necessary. 
In particular, when  $\nu_{0} > 0$, we get the following further regularity for $\theta_{n}$: 
\begin{equation}\label{6-5.5}
\left\{ \begin{array}{l}
\ds \{\theta_{n}\}_{n=1}^{\infty} \mbox{ is bounded in } L^{\infty}(0,1;H^{1}(\Omega)), \vspace{3mm}\\
\ds \theta_{n} \to \theta_{\infty} \mbox{ weakly-}\ast \mbox{ in } L^{\infty}(0,1;H^{1}(\Omega)), \vspace{3mm}\\
    \ds \alpha(\v_{n})\nabla\theta_{n} \to \alpha(\v_{\infty})\nabla\theta_{\infty}, \ts\sqrt{\beta(\v_{n})}\nabla\theta_{n} \to \sqrt{\beta(\v_{\infty})}\nabla\theta_{\infty}\vspace{3mm}\\
\ds \hspace{5mm} \mbox{ weakly in } L^{2}(0,1; L^{2}(\Omega)^{N}), 
\end{array} \right. 
\end{equation}
as $n \to \infty$. Here, we set $\{u_{n}\}_{n=1}^{\infty} := \{u(\cdot+t_{n})\}_{n=1}^{\infty}$. By (S2) and (S3), the sequence $\{ [\v_{n},\theta_{n}] \}_{n=1}^{\infty}$ satisfies 
\begin{equation}\label{6-6}
\begin{array}{l}
\ds \int_{0}^{1} ( (\v_{n})_{t}(t), \v_{n}(t)-{\bm\varpi})_{L^{2}(\Omega)^{2}} dt + \int_{0}^{1} (\nabla \v_{n}(t), \nabla (\v_{n}(t)-{\bm\varpi}))_{L^{2}(\Omega)^{2}} dt \vspace{3mm}\\
\ds + \int_{0}^{1} ([\nabla G](u_{n};\v_{n})(t), \v_{n}(t)-{\bm\varpi})_{L^{2}(\Omega)^{2}} dt + \int_{0}^{1} \int_{\Omega} \gamma(\v_{n}(t)) dxdt \vspace{3mm}\\
\ds + \int_{0}^{1}\int_{\Omega} d[(\v_{n}(t)-{\bm\varpi}) \cdot [\nabla\alpha](\v_{n}(t)) |D\theta_{n}(t)|] dt \vspace{3mm}\\
\ds + \nu_{0}^{2}\int_{0}^{1} \int_{\Omega} (\v_{n}(t) - \bm{\varpi}) \cdot [\nabla\beta](\v_{n}(t)) |\nabla \theta_{n}(t)|^{2}dxdt 
\ds 
 \le \int_{0}^{1} \int_{\Omega} \gamma({\bm\varpi}) dxdt 
\end{array}
\end{equation}
for any ${\bm \varpi} \in H^{1}(\Omega)^{2} \cap L^{\infty}(\Omega)^{2}$ and $n \in \N$, and 
\begin{equation}\label{6-7}
\begin{array}{l}
\ds \int_{0}^{1} (\alpha_{0}(\v_{n}(t)) (\theta_{n})_{t}(t), \theta_{n}(t))_{L^{2}(\Omega)} dt + \int_{0}^{1} \Phi_{\nu_{0}}(\v_{n}(t);\theta_{n}(t))dt \vspace{3mm}\\
\ds \le \int_{0}^{1} \Phi_{\nu_{0}}(\v_{n}(t);0)dt = 0
\end{array}
\end{equation}
for any $n \in \N$. Owing to (\ref{6-1.4}), (\ref{6-5}), (\ref{6-5.5}), (\ref{6-7}), and Lemma \ref{lem-gamma2}, we can see that 
\begin{equation}\label{6-8}
\begin{array}{l}
\ds 0 \le \int_{0}^{1}\int_{\Omega} d[\alpha(\v(t))|D\theta(t)|]dt + \nu_{0}^{2} \int_{0}^{1} \int_{\Omega} \beta(\v(t))|\nabla\theta(t)|^{2} dxdt \vspace{3mm}\\
\ds \hspace{3mm} \le \liminf_{n \to \infty} \int_{0}^{1} \left( \int_{\Omega} d[\alpha(\v_{n}(t))|D\theta_{n}(t)|] + \nu_{0}^{2} \int_{\Omega} \beta(\v_{n}(t)) |\nabla\theta_{n}(t)|^{2}dx \right) dt \vspace{3mm}\\
\ds \hspace{3mm} \le \limsup_{n\to\infty} \int_{0}^{1} \Phi_{\nu_{0}}(\v_{n}(t);\theta_{n}(t))dt \vspace{3mm}\\
\ds \hspace{3mm} \le - \lim_{n\to\infty} \int_{0}^{1}(\alpha_{0}(\v_{n}(t)) (\theta_{n})_{t}(t), \theta_{n}(t))_{L^{2}(\Omega)} dt = 0.
\end{array}
\end{equation}
By (A2), the above inequality means (i-c). 

Finally, using (\ref{6-5}) and (\ref{6-8}), we apply Lemma \ref{lem6-1} with $I = (0,1)$, $\rho = \alpha(\v_{\infty})$, $\{\rho_{n}\}_{n=1}^{\infty} = \{\alpha(\v_{n})\}_{n=1}^{\infty}$, $\zeta = \theta_{\infty}$, $\{\zeta_{n}\}_{n=1}^{\infty} = \{\theta_{n}\}_{n=1}^{\infty}$, $\omega = {\bm\varpi} \cdot [\nabla\alpha](\v_{\infty})$, and $\{\omega_{n}\}_{n=1}^{\infty} = \{ {\bm\varpi} \cdot [\nabla\alpha](\v_{n})\}_{n=1}^{\infty}$. Then, we see that  
\begin{equation}\label{6-9}
\begin{array}{l}
\ds \int_{0}^{1}\int_{\Omega} d[{\bm \varpi} \cdot [\nabla\alpha](\v_{n}(t))|D\theta_{n}(t)|]dt \to \int_{\Omega} d[{\bm \varpi} \cdot [\nabla\alpha](\v_{\infty})|D\theta_{\infty}|] = 0
\end{array}
\end{equation}
as $n\to \infty$, for any ${\bm\varpi} \in H^{1}(\Omega)^{2} \cap L^{\infty}(\Omega)^{2}$. 
In particular, when $\nu_{0}>0$, (\ref{6-5.5}) and (\ref{6-7}) lead to
\begin{equation}\label{6-10}
|\nabla \theta_{\infty}|_{L^{2}(\Omega)^{N}}^{2} = \lim_{n \to \infty} |\nabla\theta_{n}|_{L^{2}(\Omega)^{N}}^{2} = 0. 
\end{equation}
By (\ref{6-5}), (\ref{6-9}), and (\ref{6-10}), letting $n \to \infty$ in (\ref{6-6}) yields that 
\begin{equation*}
\begin{array}{l}
\ds (\nabla \v_{\infty}, \nabla(\v_{\infty}-{\bm\varpi}))_{L^{2}(\Omega)^{2}} + ([\nabla G](u_{\infty};\v_{\infty}), \v_{\infty}-{\bm\varpi})_{L^{2}(\Omega)^{2}} \vspace{3mm}\\
\ds \hspace{5mm} + \int_{\Omega} \gamma(\v_{\infty}) dx 
\ds \le \int_{\Omega} \gamma({\bm\varpi}) dx, 
\end{array}
\end{equation*}
for any ${\bm \varpi} \in H^{1}(\Omega)^{2} \cap L^{\infty}(\Omega)^{2}$. Hence, the proof of Main Theorem 2 is completed. 
\hfill $\Box$ 
\begin{rem}
    \begin{em}
Up to the setting of the system (S)$_{\nu}$, we can observe the convergence of the orbit $ \{ [\v(t), \eta(t)] \} = \{ [w(t), \eta(t), \theta(t)] \} $, as $ t \to \infty $, 
without taking the time sequence  $\{t_{n}\}$.  For instance, under the setting:
\begin{equation*}
\gamma(\cdot) = I_{[0,1]}(\cdot), g(\v) = \frac{1}{2}(w-\eta)^{2} - \frac{c}{2}w^{2}, \alpha(w,\eta) = \alpha(\eta), \beta(w,\eta) = \beta(w) \mbox{ and } u_{\infty} \notin [0,1],
\end{equation*}
it was shown in \cite{SW15} that the $\omega$-limit set is a singleton $ \{ \v_{\infty} \} $ of a constant vector $ \v_{\infty} = [w_\infty, \eta_\infty] \in [0, 1]^2 $, and moreover, $ [w_{\infty}, \eta_{\infty}] = [1, 1] $ (resp. $ [w_{\infty}, \eta_{\infty}] = [0, 0] $) if $ u_{\infty} < 0 $ (resp. $ u_{\infty} > 1 $).
\end{em}
\end{rem}

\end{document}